\documentclass[12pt]{article}
 \usepackage{amssymb,amsmath}
\newtheorem{theorem}{Theorem}

\newtheorem{proposition}[theorem]{Proposition}
\newtheorem{lemma}[theorem]{Lemma}
\newtheorem{corollary}[theorem]{Corollary}

\newtheorem{definition}[theorem]{Definition}

\newcommand{\R}{\mathbb{R}}
\newcommand{\Ee}{\mathbb {E}}
\newcommand{\Q}{\mathbb{Q}}

\newcommand{\Sf}{\mathbb{S}}

\newcommand{\Le}{\mathbb{L}}
\newcommand{\Ve}{\mathbb{V}}
\newcommand{\C}{\mathbb{C}}
\newcommand{\Hy}{\mathbb{H}}

\newcommand{\spa}{\mbox{span}}

\newcommand{\po}{{\hspace*{-1ex}}{\bf .  }}

\def\a{\alpha}

\def\<{\langle}
\def\n{\nabla}
\def\>{\rangle}
\def\a{\alpha}

\def\e{\epsilon}
\def\d{\partial}
\def\bea{\begin{eqnarray*} }
\def\eea{\end{eqnarray*} }
\def\be{\begin{equation} }
\def\ee{\end{equation} }

\def\proof{\noindent{\it Proof: }}
\def\qed{\ifhmode\unskip\nobreak\fi\ifmmode\ifinner
\else\hskip5 pt \fi\fi\hbox{\hskip5 pt \vrule width4 pt
height6 pt  depth1.5 pt
\hskip 1pt }}
\begin{document}

\title{ Umbilical submanifolds of $\Sf^n\times \R$.}

\author { Bruno Mendon\c ca \& Ruy Tojeiro}
\date{}
\maketitle
\begin{abstract}
We give a complete  classification of  umbilical submanifolds of arbitrary dimension and codimension of 
$\Sf^n\times \R$, extending the classification of umbilical surfaces in $\Sf^2\times \R$ by Rabah-Souam  and Toubiana  as well as  the local description of  umbilical hypersurfaces in $\Sf^n\times \R$ by Van der Veken and Vrancken.  We prove that, besides small spheres in a slice, up to isometries of the ambient space they come in a two-parameter family of rotational submanifolds 
whose  substantial codimension is  either one or two and whose profile is a  curve  in a totally geodesic $\Sf^1\times \R$ or $\Sf^2\times \R$, respectively, the former case arising in a one-parameter family. All of them are diffeomorphic to a sphere, except for a single element that is diffeomorphic to Euclidean space. We obtain explicit parametrizations of all such submanifolds. We also study more general classes of submanifolds of $\Sf^n\times \R$ and $\Hy^n\times \R$. In particular, we give a complete description of all submanifolds in those product spaces 
for which the tangent component of a unit vector field spanning the factor $\R$ is an eigenvector of all shape operators. We show that surfaces with parallel mean curvature vector in $\Sf^n\times \R$ and $\Hy^n\times \R$ having this property are rotational surfaces.
We  also prove a Dajczer-type reduction of codimension theorem for submanifolds of $\Sf^n\times \R$ and $\Hy^n\times \R$.
\end{abstract}

\section{Introduction} Roughly speaking, a submanifold of a Riemannian manifold is \emph{totally umbilical}, or simply \emph{umbilical}, if it is equally curved in all tangent directions. More precisely, an isometric immersion $f\colon\, M^m\to \tilde{M}^n$ between Riemannian manifolds  is  umbilical if there exists a normal vector field $\zeta$ along $f$ such that its second fundamental form $\alpha_f\colon\, TM\times TM\to N^fM$ with values in the normal bundle satisfies $\alpha_f(X,Y)=\<X,Y\>\zeta$ for all $X,Y\in TM$. 

 Umbilical submanifolds are the simplest submanifolds after the totally geodesic ones (for which the second fundamental form vanishes identically), and  their knowledge sheds light on the geometry of the ambient space.

Apart from space forms, however, there are few Riemannian manifolds for which  umbilical submanifolds are classified. Recently, this was accomplished  for all three-dimensional  Thurston geometries of non-constant curvature as well as for the Berger spheres in \cite{st}.
The richest case turned out to be that of the product spaces $\Sf^2\times \R$ and $\Hy^2\times \R$. For these manifolds,  it was shown that, up to isometries of the ambient space, umbilical nontotally geodesic surfaces  come in a one-parameter family of  rotational surfaces, whose  profile curves have been completely determined in terms of solutions of a certain ODE.

A local description of  umbilical hypersurfaces of $\Sf^n\times \R$ and $\Hy^n\times \R$ of any dimension $n$ was given in \cite{vv} and \cite{ckv}, respectively. Again,  the nontotally geodesic ones are  rotational hypersurfaces over curves in totally geodesic products $\Sf^1\times \R$ and $\Hy^1\times \R$, respectively.

In this paper we give a complete  classification of  umbilical submanifolds of arbitrary dimension and codimension of $\Sf^n\times \R$. 
To state our result, for a given integer $m\geq 2$ let $\phi\colon\,\Sf^{m+1}\times \R\to \R^{m+2}\setminus \{0\}$ be the conformal diffeomorphism given by $\phi(x,t)=e^tx$. Choose a closed half-line $\ell:=\{\bar x\}\times [0, \infty)\subset\R^{m+2}=\R^{m+1}\times \R$ with $\bar x\neq 0$. Let  $M^{m}_{r,h}$ be the image by $\phi^{-1}$ of the  $m$-dimensional sphere $\Sf^m_{r,h}$ in $\R^{m+2}$ of radius $r$  centered on $\ell$ that lies in the affine hyperplane through $(\bar x,h)$ orthogonal to $\ell$, with the origin removed if $h=0$ and $r=d:=|\bar x|$. Then we prove:
\begin{theorem} \po\label{thm:main2} 
The submanifold $M^m_{r,h}$ is a complete umbilical submanifold of $\Sf^{m+1}\times \R$ for every $r>0$ and $h\geq  0$. Moreover, it has the following properties:
\begin{itemize} 
   \item[$(i)$] It is diffeomorphic to  $\Sf^m$ if $(r,h)\neq (d,0)$ and to $\R^m$ if $(r,h)=(d,0)$;
   \item[$(ii)$] It lies in a totally geodesic hypersurface $\Sf^{m}\times \R\subset\Sf^{m+1}\times \R$ if and only if 
   $h=0$;  
   \item[$(iii)$] 
   $M^m_{r,0}$ is homologous to zero in $\Sf^{m}\times \R$ if $r<d$ and inhomologous to zero  if $r>d$;
   \item[$(iv)$] It is a rotational submanifold whose profile is a curve  in a totally geodesic submanifold $\Sf^2\times \R$
   (respectively,  $\Sf^1\times \R$) if $h\neq 0$ (respectively, $h=0$); 
   \item[$(v)$] $M^m_{r,h}$ is not congruent to $M^m_{r',h'}$ if $(r,h)\neq (r',h')$.
   \end{itemize}
Conversely, any  umbilical nontotally geodesic submanifold of $\Sf^n\times \R$ with dimension $m\geq 2$ is, up to an isometry of the ambient space, an open subset of  one of the following: 
   \begin{itemize}
   \item[$(i)$] a small sphere in $\Sf^n\times \{0\}$;
\item[$(ii)$]   $M^m_{r,0}$ for some $r>0$ if $n=m$;
\item[$(iii)$]  $M^m_{r,h}$ for some $r>0$ and $h\geq 0$ if $n=m+1$;
\item[$(iv)$] $M^m_{r,h}$ in a totally geodesic $\Sf^{m+1}\times \R$ for some $r>0$ and $h\geq 0$ if $n> m+1$.
\end{itemize}
  \end{theorem}
  
  Moreover, we provide an explicit parametrization of all submanifolds $M^m_{r,h}$, $r>0$, $h\geq 0$ (see Proposition \ref{prop:paramet} below) in terms of elementary functions. The precise meaning of $M^m_{r,h}$ being \emph{rotational} is explained in Section $4$.

   In the process of proving Theorem \ref{thm:main2}, we have been led to study more general classes of submanifolds with interest on their own. 

Let 
 $\Q_\e^n$ denote either  $\Sf^n$, $\R^n$ or  $\Hy^n$, according as $\e=1$, $\e=0$ or $\e=-1$, respectively. Given an isometric immersion  $f\colon\,M^m\to \Q_\e^n\times\R$,  
 let   
 $\frac{\partial}{\partial t}$ be a unit
vector field tangent to the second factor. Thus, for $\epsilon=0$ we just choose a unit constant vector field $\frac{\partial}{\partial t}$ in $\R^{n+1}$.
Then, a tangent vector field  $T$ on $M^m$ and a normal vector field  $\eta$ along $f$ are defined by
\be\label{eq:ddt} \frac{\partial}{\partial t}={f}_* T+ \eta.\ee

 We denote by ${\cal A}$ the class of  isometric immersions $f\colon\,M^m\to \Q_\e^n\times\R$ with the property that $T$ is an eigenvector of all shape operators of $f$. Our  next result is a complete description of all isometric immersions in class ${\cal A}$.  
First note that trivial examples are products $N^{m-1}\times \R$, where $N^{m-1}$ is a submanifold of $\Q_\e^n$, which correspond to the case in which the normal vector field $\eta$ in (\ref{eq:ddt}) vanishes identically. We call these examples \emph{vertical cylinders}. More interesting ones are constructed  as follows. We consider the case $\epsilon\in \{-1, 1\}$, the case $\epsilon=0$ being similar.

Let   $g\colon\, N^{m-1}\to \Q_\e^n$  be an isometric immersion. Assume that there exists an orthonormal  set of  parallel  normal vector fields $\xi_1, \ldots, \xi_k$ along $g$. This assumption is  satisfied, for instance,  if $g$ has flat normal bundle.  Thus, the  vector subbundle $E$ with rank $k$ of the normal bundle $N^gN$ of $g$ spanned by $\xi_1, \ldots, \xi_k$ is  parallel and  flat. 
Let $j\colon\,\Q_\e^n\to \Q_\e^n\times \R$ and $i\colon\,\Q_\e^n\times\R\to  \Ee^{n+2}$ denote the canonical inclusions, and let $k=i\circ j$.  Here $\Ee^{n+2}$ denotes either Euclidean space $\R^{n+2}$ if $\e=1$ or Lorentzian space $\Le^{n+2}$ if $\e=-1$. Set $\tilde{\xi_i}=k_*\xi_i$, $1\leq i\leq k$,  $\tilde{\xi}_0=\tilde{g}:=k\circ g$ and  $\tilde{\xi}_{k+1}=i_*\d/\d t$. Then the vector subbundle $\tilde E$ of $N^{\tilde g}N$ whose fiber $\tilde{E}(x)$ at $x\in N^{m-1}$ is spanned by $\tilde\xi_0, \ldots, \tilde \xi_{k+1}$ is also parallel and flat,  and we may define a vector bundle isometry
$\phi\colon\, N^{m-1}\times \Ee^{k+2}\to \tilde E$ by
$$\phi_x(y):=\phi(x,y)=\sum_{i=0}^{k+1}y_i\tilde \xi_i, \,\,\,\,\mbox{for}\,\,\,y=(y_0, \ldots, y_{k+1})\in \Ee^{k+2}.$$  
Now let 
$$f\colon\, M^m:=N^{m-1}\times I\to \Q_\e^n\times \R$$  be given by 
\be\label{eq:f}\tilde{f}(x,s):=(i\circ f)(x,s)=\phi_x(\alpha(s))=\sum_{i=0}^{k+1}\a_i(s)\tilde{\xi}_i(x),\ee
 where 
 $\alpha\colon\, I\to \Q_\e^k\times \R\subset \Ee^{k+2}$, $\alpha=(\a_0,\ldots, \a_k, \a_{k+1})$,   is a smooth regular curve such that  $\e \a_0^2+\a_1^2+\ldots+ \a_k^2=\e$ and $\a_{k+1}$ has nonvanishing derivative. Notice that vertical cylinders correspond to the case in which the curve $\alpha$ is the generator of $\Q_\e^k\times \R$ through $(1,0,\cdots,0)\in \Q_\e^k$.

\begin{theorem}\po\label{thm:main} The map $f$ defines, at  regular points,   an immersion 
in class ${\cal A}$.
Conversely,  any isometric immersion  $f\colon\, M^m\to \Q_\e^n\times \R$, $m\geq 2$, 
in class ${\cal A}$ 
is locally given in this way.
\end{theorem}

A necessary and sufficient condition for a point $(x,s)\in M^m=N^{m-1}\times I$ to be regular for $f$ is given in part $(ii)$ of 
Proposition \ref{prop:cu} below.

The map $\tilde{f}$ is a \emph{partial tube over $\tilde g$ with type fiber $\alpha$} in the sense of \cite{cw} (see also \cite{cu}). 
Geometrically, $\tilde{f}(M)$ is obtained by parallel transporting the curve $\alpha$ in a product submanifold $\Q_\epsilon^k\times \R$ of a fixed normal  space of $\tilde{g}$ with respect to its normal connection.

Theorem \ref{thm:main} extends to submanifolds with arbitrary codimension the main result of \cite{to}, where the case of  hypersurfaces was studied. That the preceding construction  coincides with the one in Theorem $1$ of \cite{to} in the hypersurface case was already observed in  Remarks $7$-$(ii)$ in that paper. Some important classes of hypersurfaces of $\Q_\e^n\times \R$, $\epsilon \in \{-1, 1\}$, that are included in class ${\cal A}$ are  hypersurfaces with constant sectional curvature \cite{mt}, rotational hypersurfaces \cite{dfv} and constant angle hypersurfaces (see, e.g., \cite{to}; see also Corollary \ref{cor:constantangle} below and the comments before it). 

Let $f\colon\, M^m\to \Q_\e^n\times \R$,  $\epsilon \in \{-1, 1\}$, be an isometric immersion, and set  $\tilde{f}=i\circ f$, where $i\colon\,\Q_\e^n\times\R\to  \Ee^{n+2}$ is the canonical inclusion. It was shown in \cite{mt} that if  $m=n$ then $f$ is in class ${\cal A}$ if and only if the vector field $T$ in (\ref{eq:ddt}) is nowhere vanishing and $\tilde{f}$ has flat normal bundle. For submanifolds of higher codimension we have the following.

\begin{corollary}\po\label{cor:fnb} 
The following assertions are equivalent:
\begin{itemize}
\item[(i)] The vector field $T$ in (\ref{eq:ddt}) is nowhere vanishing and $\tilde f$ has flat normal bundle;
\item[(ii)] $f$ has flat normal bundle and is in class ${\cal A}$; 
\item[(iii)] $\tilde f$ is locally given as in (\ref{eq:f}) in terms of an isometric immersion $g\colon\, N^{m-1}\to \Q_\e^n$ with flat normal bundle and a smooth regular curve $\alpha\colon\, I\to \Q_\e^k\times \R\subset \Ee^{k+2}$, $\alpha=(\a_0,\ldots, \a_k, \a_{k+1})$,   with  $\a'_{k+1}$ nowhere  vanishing.
\end{itemize}
\end{corollary}

Observe that  the vector field  $T$ vanishes at some point if and only if $f(M^m)$ is tangent to the slice $\Q_\e^n\times \{t\}$ of  $\Q_\e^n\times \R$ through that point. If $T$ vanishes on an open subset $U\subset M^m$ then $f(U)$ is contained in some slice. 

Notice that a surface in  class ${\cal A}$ has automatically flat normal bundle. Hence, by Corollary \ref{cor:fnb}, a surface in
$\Q_\e^n\times \R$, $\epsilon \in \{-1, 1\}$,  is in  class ${\cal A}$ if and only if it has flat normal bundle as a surface in the underlying flat space $\Ee^{n+2}$ (and is nowhere tangent to a slice). By Theorem \ref{thm:main}, any such surface is given by (\ref{eq:f}) in terms of a unit-speed curve  $g\colon\, J\to \Q_\e^n$ and a smooth regular curve $\alpha\colon\, I\to \Q_\e^k\times \R\subset \Ee^{k+2}$, $\alpha=(\a_0,\ldots, \a_k, \a_{k+1})$,   with  $\a'_{k+1}$ nowhere  vanishing. Clearly,  in this case the existence of an orthonormal set of parallel normal vector fields $\xi_1, \ldots , \xi_k$ along $g$ is automatic for any $1\leq k\leq n-1$.

 In the case of a hypersurface $f\colon\,M^{n}\to \Q_\e^n\times \R$, the vector field $\eta$ in  (\ref{eq:ddt}) can be written as $\eta=\rho N$, where $N$ is a unit normal vector field along $f$. Then $f$ is called a {\em constant angle hypersurface\/} if the function $\rho$ is constant on $M^n$. One possible  way to generalize this notion to submanifolds of higher codimension is to require the vector field $\eta$ to be parallel in the normal connection. It turns out that submanifolds with this property also belong to class ${\cal A}$, and this leads to the following classification of them,  extending Corollary $2$ in \cite{to}.
  
\begin{corollary}\po\label{cor:constantangle} Let $f\colon\, M^m:=N^{m-1}\times I\to \Q_\e^n\times \R$ be given  by
(\ref{eq:f}) with $\alpha\colon\, I\to \Q_\e^k\times \R$  a geodesic of $\Q_\e^k\times \R$. Then $f$ defines, at  regular points,   an immersion for which the vector field $\eta$ in (\ref{eq:ddt}) is parallel in the normal connection. 
Conversely,  any isometric immersion   $f\colon\, M^m\to \Q_\e^n\times \R$, $m\geq 2$, such that $T$ is nowhere vanishing and $\eta$  is parallel in the normal connection 
is locally given in this way.
\end{corollary}

Another important subclass of class ${\cal A}$ is that of \emph{rotational submanifolds} in $\Q_\epsilon^n\times \R$ with curves in totally geodesic submanifolds $\Q_\epsilon^\ell\times \R\subset \Q_\epsilon^n\times \R$ as profiles (see Section~$4$).  We obtain the following characterization of independent interest of them.
 
 \begin{corollary}\po\label{prop:rot} Let $f\colon\,M^{m}\to \Q_\e^n\times \R$, $\epsilon\in \{-1, 1\}$,  be an isometric immersion. Then the following assertions are equivalent:
 \begin{itemize}
 \item[$(i)$] $f$ is a rotational submanifold whose profile is a curve in a totally geodesic submanifold $\Q_\e^{n-m+1}\times\R\subset \Q_\epsilon^n\times \R$;
 \item[$(ii)$] $f$ is given  as in (\ref{eq:f})  in terms of an umbilical isometric immersion $g\colon\,N^{m-1}\to \Q_\epsilon^n$ (a geodesic circle, if $m=2$);
 \item[$(iii)$] there exists a normal vector field $\zeta$ along $f$ such that 
 \be\label{eq:sff0}\alpha_f(X,Y)=\<X,Y\>\zeta\,\,\,\mbox{for all}\,\,\,X\in TM\,\,\mbox{and}\,\,Y\in \{T\}^\perp,\ee
 where $T$ is the vector field defined by (\ref{eq:ddt}), and $\zeta$ is parallel in the normal connection along $\{T\}^\perp$ if $m=2$.
 \end{itemize}
 Moreover, if $\e=1$ then the preceding assertions are equivalent to $f$ being given as in (\ref{eq:f})  in terms of a totally geodesic  isometric immersion $g\colon\,N^{m-1}\to \Q_\epsilon^n$. 
 This is also the case if $\e=-1$ and $f$ is assumed to be of hyperbolic type in $(i)$ and $g$ an equidistant hypersurface  in $(ii)$.
 \end{corollary}

     Notice that in the hypersurface case, i.e., for $n=m$,   the second fundamental form satisfies (\ref{eq:sff0})  if and only if  $f$ has at most two principal curvatures, and if it has exactly two then one of them is simple with $T$ as an eigenvector.

 A key step in the classification of umbilical submanifolds in $\Sf^n\times \R$ is the following result on reduction of codimension of isometric immersions into $\Q_\e^n\times \R$. That an  isometric immersion $f\colon\, M^m\to \Q_\e^n\times \R$ \emph{reduces codimension to} $p$, or has \emph{substantial codimension} $p$,  means that $f(M^n)$ is contained in a totally geodesic submanifold $\Q_\e^{m+p-1}\times \R$ of $\Q_\e^n\times \R$. 
We denote by $\nabla^\perp$ the normal connection of $f$ and by $N_1(x)$ the \emph{first normal space} of $f$ at $x$, i.e., the subspace of 
$N^f_xM$  spanned by its second fundamental form.
 
\begin{lemma}\po\label{le:redcod} Let $f\colon\,M^m\to \Q_\e^n\times \R$, $\epsilon\in \{-1, 1\}$, be an isometric immersion. Let $\eta$ be the normal vector field defined by (\ref{eq:ddt}). Assume that $L:=N_1 + \spa\{\eta\}$ is a subbundle of $N^fM$ with rank $\ell<n+1-m$ and that  $\nabla^\perp N_1\subset L$. Then $f$ reduces codimension to~$\ell$. 
 \end{lemma}
 
 Lemma \ref{le:redcod} should be compared with its well-known counterpart for submanifolds of space forms (see, e.g., \cite{er}), in which case the corresponding condition for a submanifold to reduce codimension is that its first normal spaces form a parallel subbundle of the normal bundle. A necessary and sufficient condition for parallelism of the first normal bundle of a submanifold of a space form in terms of  its normal curvature tensor $R^\perp$ and mean curvature vector field $H$  was obtained by Dajczer \cite{da} (see also Chapter $4$ of \cite{daa}). The proof of Dajczer's theorem can be easily adapted to yield the following  result for submanifolds of  $\Q_\e^n\times \R$.

 \begin{theorem}\po\label{thm:dajczer} Let $f\colon\,M^m\to \Q_\e^n\times \R$, $\epsilon\in \{-1, 1\}$, be an isometric immersion. Let $\eta$ be the normal vector field defined by (\ref{eq:ddt}). Assume that $L:=N_1 + \spa\{\eta\}$ is a subbundle of $N^fM$ of rank $\ell<n+1-m$. Then  $\nabla^\perp N_1\subset L$ if and only if  the following two conditions hold:
 \begin{itemize}
 \item[(i)] $\nabla^\perp R^\perp|_{L^\perp}=0$;
 \item[(ii)] $\nabla^\perp H\in L$.
 \end{itemize}
 \end{theorem}
 
 As  an  application of Theorem \ref{thm:dajczer}, in Subsection $5\!-\!1$ we give a simple proof of  Theorem~$1$  in \cite{adct} on surfaces with parallel mean curvature vector in $\Q_\epsilon^n\times \R$. By using this result together with Corollary \ref{prop:rot} we prove the following.

     \begin{corollary}\po \label{cor:pmc} Any surface $f\colon\,M^{2}\to \Q_\e^n\times \R$, $\epsilon \in \{-1, 1\}$,  in class ${\cal A}$ with parallel mean curvature vector is a rotational surface in a totally geodesic submanifold $\Q_\epsilon^m\times \R$, $m\leq 4$, over a curve in a totally geodesic submanifold $\Q_\epsilon^s\times \R$, $s\leq 3$.
     \end{corollary}
     
     In the case $n=2$, the preceding corollary is a special case of Theorem $3$ in \cite{to}, in which hypersurfaces $f\colon\,M^{n}\to \Q_\e^n\times \R$ in class ${\cal A}$ with constant mean curvature and arbitrary dimension $n$ were classified. That for $n=2$ they are all rotational surfaces was pointed out in  part $(i)$ of Remarks $7$ in that paper. Recently we learned that this was independently proved in Theorem $1$ of  \cite{adct2}.
     
     In \cite{adct0}, the authors introduced the real quadratic form 
$$Q(X, Y ) = 2\<\alpha(X, Y ),H\> - \epsilon\<X,T\>\<Y,T\>$$
on a surface $f\colon\,M^{2}\to \Q_\e^n\times \R$, 
as a generalization to higher codimensions of the Abresch--Rosenberg real quadratic form defined in \cite{ar}. 
Extending the result in \cite{ar} for constant mean curvature surfaces, they proved that 
the $(2,0)$-part $Q^{(2,0)}$ of $Q$ is holomorphic for surfaces with parallel mean curvature vector field. This means that
if $(u,v)$ are isothermal coordinates on $M^2$, then the complex function 
$$Q(Z,Z)=2\<\alpha(Z, Z),H\> - \epsilon\<Z,T\>^2$$
is holomorphic, where $Z=\frac{1}{\sqrt{2}}(\frac{\d}{\d u}+i\frac{\d}{\d v})$ and the metric on $M^2$ is extended to a $\C$-bilinear map.

The same authors  observed in \cite{adct} that surfaces with parallel mean curvature vector that are not contained in a slice of $\Q_\e^n\times\R$ and for which $Q^{(2,0)}$ vanishes identically belong to  class ${\cal A}$. They also proved that a surface 
$f\colon\,M^2\to \Q_\e^n$ with parallel mean curvature vector field has vanishing $Q^{(2,0)}$ if either $M^2$ is homeomorphic to a sphere 
or if $\epsilon=-1$, $K\geq 0$ and  $K$ is not identically zero. By means of  Corollary \ref{cor:pmc}, we obtain the following improvement of  the conclusions in part $4$ of both
Theorems $2$ and $3$ of \cite{adct}:

\begin{corollary}\po \label{cor:pmc2} Let $f\colon\,M^{2}\to \Q_\e^n\times \R$, $\epsilon \in \{-1, 1\}$,  be a surface  with parallel mean curvature vector. Suppose  $f(M^2)$ is not contained in a slice of $\Q_\e^n\times \R$ and  either
\begin{itemize}
\item[$(i)$] $M^2$ is homeomorphic to a sphere, or
\item[$(ii)$] $\epsilon=-1$, $M^2$ is complete with $K\geq 0$ and  $K$ is not identically zero.
\end{itemize}
Then $f$ is a rotational surface (of spherical type in case $(i)$) in a totally geodesic submanifold $\Q_\epsilon^m\times \R$, $m\leq 4$, over a curve in a totally geodesic submanifold $\Q_\epsilon^s\times \R$, $s\leq 3$.
\end{corollary}

We observe that, if $f(M^2)$ is  contained in a slice $\Q_\epsilon^n\times \{t\}$ of $\Q_\e^n\times \R$, then by Theorem~$4$ in \cite{yau} either  $f$ is a minimal surface of an  umbilical hypersurface of $\Q_\epsilon^n\times \{t\}$ or it is a  surface with constant mean curvature in a three-dimensional  umbilical or totally geodesic submanifold of $\Q_\epsilon^n\times \{t\}$. Moreover, if 
$M^2$ is homeomorphic to a sphere, then by Hopf's Theorem it must be a totally umbilical $2$-sphere of  $\Q_\e^n\times \{t\}$.
 
 The paper is organized as follows. In the next section we recall the basic equations of an isometric immersion into $\Q_\e^n\times \R$. In Section $3$ we study submanifolds in class ${\cal A}$ and prove Theorem \ref{thm:main} as well as Corollaries \ref{cor:fnb} and \ref{cor:constantangle}. Section $4$ is devoted to rotational submanifolds. In particular, Corollary \ref{prop:rot} is proved. In Section $5$ we prove Lemma~\ref{le:redcod} and Theorem \ref{thm:dajczer} on reduction of codimension of isometric immersions into $\Q_\e^n\times \R$. Then we apply the latter to  give  a simple proof of  Theorem $1$  in \cite{adct} on surfaces with parallel mean curvature vector in $\Q_\epsilon^n\times \R$.  We conclude this section with the proof of Corollary~\ref{cor:pmc}.  Finally, in the last section we prove  Theorem \ref{thm:main2} on the classification of umbilical submanifolds of $\Sf^n\times \R$.

\section{Preliminaries}

      In this  section we recall the fundamental equations of an isometric immersion $f\colon\,M^{m}\to \Q_\e^n\times \R$. 
      
        Using that
$\frac{\partial}{\partial t}$ is a parallel vector field in $\Q_\e^n\times\R$, we obtain by differentiating (\ref{eq:ddt}) that 
\begin{eqnarray}\label{eq:NablaT}
\nabla_XT= A^f_\eta X
\end{eqnarray}
and 
\begin{eqnarray}\label{eq:Derivadaeta}
\alpha_f(X,T)=-\nabla_X^\perp \eta,
\end{eqnarray}
for all $X\in TM$. Here and in the sequel $A_\eta^f$ stands for the shape operator of $f$ in the direction $\eta$, given by
$$\<A_\eta^fX, Y\>=\<\alpha_f(X,Y), \eta\>\,\,\,\,\mbox{for all}\,\,\,X,Y\in TM.$$

Notice  that the vector field $T$ is a gradient vector field. Namely, if $\epsilon\in \{-1, 1\}$ and $\tilde{f}=i\circ f$, where $i\colon\,\Q_\e^n\times\R\to  \Ee^{n+2}$ denotes the canonical inclusion, then $T$  is  the gradient of the height function $h=\<\tilde{f},i_*\frac{\partial}{\partial t}\>$. If $\epsilon=0$ then $T$ is the gradient of $h=\<f,\frac{\partial}{\partial t}\>$.

 The Gauss, Codazzi and Ricci equations for $f$ are, respectively (see, e.g., \cite{ltv}),  
 \be\label{eq:gauss2}
 R(X,Y)W= \epsilon(X\wedge Y-\<Y,T\>X\wedge T+\<X,T\>Y\wedge T)W +A^f_{\alpha(Y,W)}X -A^f_{\alpha(X,W)}Y,\ee
\be\label{eq:codazzi2}
\left(\n^\perp_X\alpha\right)(Y,W)-\left(\n^\perp_Y\alpha\right)(X,W)=\epsilon (\<X,W\>\<Y,T\>-\<Y,W\>\<X,T\>)\eta
\ee
and 
\be\label{eq:ricci2}
 R^\perp(X,Y)\zeta=\a(X,A^f_\zeta Y)-\a(A^f_\zeta X,Y).\ee
Equation (\ref{eq:codazzi2}) can also be written as
\be\label{eq:codazzi3}
(\nabla_Y A^f)(X, \zeta)- (\nabla_X A^f)(Y, \zeta)=\epsilon\<\eta, \zeta\>(X\wedge Y)T,
\ee
where $(X\wedge Y)T=\<Y,T\>X-\<X,T\>Y$.

Although this will not be used in the sequel, it is worth mentioning that equations (\ref{eq:NablaT}) --(\ref{eq:ricci2}) completely determine an isometric immersion $f\colon\,M^{m}\to \Q_\e^n\times \R$ up to isometries of $\Q_\e^n\times \R$ (see Corollary $3$ of \cite{ltv}).

We now relate the second fundamental forms and normal connections of $f$ and $\tilde{f}$. First notice that $\hat \nu=\pi\circ i$ is a unit normal vector field to the inclusion $i\colon\,\Q_\e^n\times\R\to  \Ee^{n+2}$, $\e\in \{-1,1\}$, where $\pi\colon\,  \Ee^{n+1}\times\R\to \Ee^{n+1}$ is the projection, and
\begin{eqnarray*}\tilde \n_Z\hat\nu&=&\pi_*i_*Z=i_*Z-\<i_*Z,i_*\frac{\d}{\d t}\>i_*\frac{\d}{\d t}\\
&=&i_*(Z-\<Z,\frac{\d}{\d t}\>\frac{\d}{\d t}),
\end{eqnarray*}
for every $Z\in T(\Q_\e^n\times\R)$, where $\tilde\n$ is the derivative in $\Ee^{n+2}$. Hence
\be\label{eq:sffi} A^i_{\hat \nu}Z=-Z+\<Z,\frac{\d}{\d t}\>\frac{\d}{\d t}.\ee

The normal spaces of $f$ and $\tilde f$ are related by 
$$N^{\tilde f} M=i_*N^fM\oplus \spa\{\nu\},$$
where $\nu=\hat \nu\circ f=\pi\circ \tilde f$. Let $\bar \nabla$ denote the Levi-Civita connection of $\Q_\e^n\times\R$. 
Given $\xi\in N^fM$, we obtain from (\ref{eq:sffi}) that 
\begin{eqnarray*} \tilde \nabla_X i_*\xi&=&i_*\bar \nabla_X\xi+\alpha_i(f_*X, \xi)\\
&=&-\tilde f_*A^f_\xi X+i_*\nabla_X^\perp\xi+\<X,T\>\<\xi, \eta\>\nu,
\end{eqnarray*}
hence
$$A^{\tilde f}_{i_*\xi}=A^f_\xi$$
 and 
\be\label{eq:nconns} \tilde \n^\perp_Xi_*\xi=i_*\n^\perp_X\xi+\<X,T\>\<\xi,\eta\>\nu\ee
for every $\xi\in N^fM$,  where $\tilde \n^\perp$ is the normal connection of $\tilde f$. On the other hand, 
$$\tilde \nabla_X\nu=\tilde \nabla_X\hat\nu\circ f=\tilde \nabla_{f_*X}\hat\nu=\tilde f_*(X-\<X,T\>T)-\<X,T\>i_*\eta,$$
hence
$$ A^{\tilde f}_\nu X=-X+\<X,T\>T,$$
or equivalently,
\be\label{eq:eigenvectorAxi}
{A}^{\tilde f}_\nu T=-\|\eta\|^2T\,\,\,\,\mbox{and}\,\,\,\,{A}^{\tilde f}_\nu X=-X\,\,\,\mbox{for}\,\,X\in \{T\}^\perp,
\ee
and
\be\label{eq:normalnu}\tilde \n^\perp_X\nu=-\<X,T\>i_*\eta.\ee

\section{Class ${\cal A}$} 

In this section  we study submanifolds in class ${\cal A}$. In particular, we give the proofs of  Theorem \ref{thm:main} and of Corollaries \ref{cor:fnb} and \ref{cor:constantangle}. We start with the following observation.

\begin{proposition}\po\label{prop:xi} Assume that the vector field $T$ in (\ref{eq:ddt}) is nowhere vanishing. Then the following assertions are equivalent:
\begin{itemize}
\item[$(i)$] $T$ is an eigenvector of $A^f_\zeta$ for all $\zeta\in N^fM$;
\item[$(ii)$] $\eta$ is parallel along $\{T\}^\perp$;
\item[$(iii)$] ${A}^{\tilde f}_\nu$ commutes with $A^f_\zeta$ for all $\zeta\in N^fM$.
\end{itemize}
\end{proposition}
\proof The equivalence between  $(i)$ and $(ii)$  follows from (\ref{eq:Derivadaeta}), whereas  (\ref{eq:eigenvectorAxi}) implies the equivalence between $(i)$ and $(iii)$.\vspace{2ex}\qed

Before going into the proof of Theorem \ref{thm:main}, we write down in the next proposition  the differential, the normal space and the second fundamental form of an immersion  $$\tilde{f}=i\circ f\colon\, M^m:=N^{m-1}\times I\to \Q_\e^n\times \R\subset \Ee^{n+2}, \,\,\,\epsilon\in \{-1, 1\},$$ which is given by (\ref{eq:f})  in terms of an isometric immersion $g\colon\, N^{m-1}\to \Q_\e^n$ and a smooth regular curve $\alpha\colon\, I\to \Q_\e^k\times \R\subset \Ee^{k+2}$, $\alpha=(\alpha_0,\ldots, \alpha_k, \alpha_{k+1})$,    with $\e \alpha_0^2+\alpha_1^2+\ldots+ \alpha_k^2=\e$. The case $\epsilon=0$ is similar. 
We use the notations before the statement of Theorem \ref{thm:main}. Given $x\in N^{m-1}$,   $X\in T_xN$ and $s\in I$, we denote by $X^{\cal H}$ the unique vector in $T_{(x,s)}M$ such that ${\pi_1}_*X^{\cal H}=X$ and ${\pi_2}_*X^{\cal H}=0$, where $\pi_1\colon\, M^m\to N^{m-1}$ and $\pi_2\colon\, M^m\to I$ are the canonical projections.

  \begin{proposition}\po \label{prop:cu} 
  The following holds:
  \begin{itemize}
  \item[(i)] The differential of $\tilde f$ is given by 
   \be\label{eq:diff1}{\tilde f}_*(x,s) X^{\cal H}=\tilde{g}_*(x)(\alpha_0(s)I-\sum_{i=1}^k\alpha_i(s)A^g_{\xi_i}(x))X,\,\,\,\mbox{for every}\,\,\, X\in T_xN,\ee
  where $I$ is the identity endomorphism of $T_xN$, 
  and \be\label{eq:diff2}{\tilde f}_*(x,s)\frac{\d}{\d s}=\phi_x(\alpha'(s)).
  \ee
  \item[(ii)] The map $\tilde{f}$ (and hence $f$)  is an immersion at $(x,s)$ if and only if 
  $$P_s(x):=\alpha_0(s)I-\sum_{i=1}^k\alpha_i(s)A^g_{\xi_i}(x)=-A^{\tilde g}_{\phi_x(\bar \alpha(s))},$$
  where $\bar{\alpha}(s)=(\alpha_0(s),\ldots, \alpha_k(s), 0)$, is an invertible endomorphism of $T_xN$.
  \item[(iii)] If $\tilde{f}$ is an immersion at $(x,s)$ then
  $$N_{(x,s)}^{\tilde{f}}{M}=k_*E(x)^\perp\oplus \phi_x(\alpha'(s)^\perp)\subset N_x^{\tilde g}N,$$
 where $E(x)^\perp$ is the orthogonal complement of $E(x)$ in $N^g_xN$,
  and \be\label{eq:N} N_{(x,s)}^{\tilde{f}}M=i_*N_{(x,s)}^{{f}}M\oplus \spa\{(\pi\circ\tilde{f})(x,s)\}=i_*N_{(x,s)}^{f}{M}\oplus \phi_x(\bar{\alpha}(s)).\ee 
  \item[(iv)] If $\tilde{f}$ is an immersion at $(x,s)$ then 
  \be\label{eq:alphapt}A^{\tilde f}_\xi(x,s) X^{\cal H}=(P_s(x)^{-1}A^{\tilde g}_\xi(x) X)^{\cal H}\ee
  for all $\xi\in N_{(x,s)}^{\tilde f}M$ and $X\in T_xN$,
  \be\label{eq:alphapt2}A^{\tilde f}_{\xi}(x,s)\frac{\d}{\d s}=0, \,\,\,\,\mbox{if}\,\,\,\,\xi \in k_*E(x)^\perp,\ee
  and 
  \be\label{eq:alphapt3}A^{\tilde f}_{\phi_x(\zeta)}(x,s)\frac{\d}{\d s}=\frac{\<\alpha''(s),\zeta\>}{\<\alpha'(s), \alpha'(s)\>}\frac{\d}{\d s},\,\,\,\mbox{if}\,\,\,\,\,\zeta\in \Ee^{k+2},\,\, \<\zeta, \alpha'(s)\>=0.\ee
Moreover, \be\label{eq:sffs}A_\zeta^f(x,s)=A_{i_*\zeta}^{\tilde f}(x,s)\ee for every $\zeta \in N_{(x,s)}^{{f}}{M}$.
  \end{itemize}
  \end{proposition} 
\proof Given a smooth curve $\gamma\colon\,J\to N^{m-1}$ with $0\in J$, $\gamma(0)=x$ and $\gamma'(0)=X$, for each $s\in I$ 
let $\gamma_s\colon\, J\to M^m$ be given by $\gamma_s(t)=(\gamma(t), s)$. Then $\gamma_s(0)=(x,s)$ and $\gamma'_s(0)=X^{{\cal H}}$.
Hence
$$\begin{array}{l}{\tilde f}_*(x,s) X^{\cal H}=\frac{d}{dt}|_{t=0}\tilde f(\gamma_s(t))=\frac{d}{dt}|_{t=0}\sum_{i=0}^{k+1}\alpha_i(s)\tilde{\xi}_i(\gamma(t))\vspace{1ex}\\\hspace*{11.6ex}=
\tilde{g}_*(x)(\alpha_0(s)I-\sum_{i=1}^k\alpha_i(s)A^{\tilde g}_{\tilde{\xi}_i}(x))X,\end{array}$$
and (\ref{eq:diff1}) follows from the fact that $A^{\tilde g}_{\tilde{\xi}_i}=A^g_{\xi_i}$ for any $1\leq i\leq k$. 

The proof of (\ref{eq:diff2}) is straightforward, and the assertions in $(ii)$ and $(iii)$ follow immediately from $(i)$.
To prove (\ref{eq:alphapt}), given $\xi\in N_{(x,s)}^{\tilde f}M$ and $X\in T_xN$, let $\gamma\colon\,J\to N^{m-1}$ and $\gamma_s\colon\, J\to M^m$ be as in the beginning of the proof. Then, using (\ref{eq:diff1}) we obtain
$$\begin{array}{l}-{\tilde f}_*(x,s)A^{\tilde f}_\xi(x,s) X^{\cal H}=(\tilde\nabla_{X^{\cal H}}\xi)^{T}=(\frac{d}{dt}|_{t=0}\xi(\gamma_s(t)))^T=
-\tilde{g}_*(x)A^{\tilde g}_\xi(x) X\vspace{1ex}\\\hspace*{21.2ex}=-\tilde{g}_*(x)P_s(x)P_s(x)^{-1}A^{\tilde g}_\xi(x) X=-{\tilde f}_*(x,s)(P_s(x)^{-1}A^{\tilde g}_\xi(x) X)^{\cal H},\end{array}$$
and (\ref{eq:alphapt}) follows. Here, putting $T$ as a superscript of a vector  means taking its  tangent component.

Formula (\ref{eq:alphapt2}) is clear. As for (\ref{eq:alphapt3}), given $\zeta\in \Ee^{k+2}$ with $\<\zeta, \alpha'(s)\>=0$, extend
$\zeta$ to a parallel normal vector field along $\alpha$, so that 
$$\zeta'(s)=\frac{\<\zeta'(s),\alpha'(s)\>}{\<\alpha'(s),\alpha'(s)\>}\alpha'(s)= -\frac{\<\alpha''(s),\zeta(s)\>}{\<\alpha'(s),\alpha'(s)\>}\alpha'(s).$$
 Then we have
$$  \begin{array}{l}\displaystyle{-{\tilde f}_*(x,s)A^{\tilde f}_{\phi_x(\zeta)}(x,s)\frac{\d}{\d s}=\tilde\nabla_{\frac{\d}{\d s}}\phi_x(\zeta)=\phi_x(\zeta'(s))=
-\frac{\<\alpha''(s),\zeta(s)\>}{\<\alpha'(s),\alpha'(s)\>}\phi_x(\alpha'(s))}\vspace{1ex}\\\hspace*{23.5ex}\displaystyle{=-{\tilde f}_*(x,s)\frac{\<\alpha''(s),\zeta(s)\>}{\<\alpha'(s),\alpha'(s)\>}\frac{\d}{\d s},}\end{array}$$
where we have used (\ref{eq:diff2}) in the last equality. This gives (\ref{eq:alphapt3}) and completes the proof, for  (\ref{eq:sffs}) is clear.\vspace{2ex}\qed

\noindent \emph{Proof of Theorem \ref{thm:main}:}  It follows from (\ref{eq:diff1}) and (\ref{eq:diff2}) that $\<X^{\cal H}, \frac{\d}{\d s}\>=0$ for any $X\in TN$, with respect to the metric induced by $f$. On the other hand, we also have 
 from (\ref{eq:diff1}) that $\<X^{\cal H},T\>=\<\tilde{f}_*X^{\cal H}, i_*\frac{\d}{\d t}\>=0$ for any $X\in TN$. Hence 
$T$ is in the direction of $\d/\d s$.
We have $$\<T,\d/\d s\>=\<\tilde f_*T, \tilde f_*\d/\d s\>=\<i_*\d/\d t, \phi_x(\alpha'(s))\>=\alpha_{k+1}'(s),$$
hence
$$T=\frac{\<T,\d/\d s\>}{\<\d/\d s, \d/\d s\>}\,\frac{\d}{\d s}=\frac{\alpha_{k+1}'(s)}{\|\alpha'(s)\|^2}\,\frac{\d}{\d s}.$$
In particular, $T$ is nowhere vanishing by the assumption that $\alpha_{k+1}'(s)\neq 0$ for all $s\in I$.  That $f$ belongs to class ${\cal A}$ now follows from (\ref{eq:alphapt2}),  (\ref{eq:alphapt3}) and (\ref{eq:sffs}).

 Let us prove the converse. Since $f\colon\, M^m\to \Q_\e^n\times \R$ belongs to class ${\cal A}$,  the vector field $T$ is nowhere vanishing, and using  (\ref{eq:NablaT}) and the fact that $T$ is a gradient vector field  we obtain
 \be\label{eq:ntt}\<\nabla_TT,X\>=\<\nabla_XT,T\>=\<A^f_\eta T,X\>=0\,\,\,\,\mbox{for any}\,\,\,X\in \{T\}^\perp.\ee Hence, the one-dimensional  distribution spanned by $T$ is totally geodesic.
  Moreover,  since $T$ is a gradient, then  the orthogonal distribution $\{T\}^\perp$ is integrable. Therefore, there exists locally a diffeomorphism $\psi\colon\,N^{m-1}\times I\to M^m$,  where $I$ is an open interval containing $0$, such that $\psi(x,\cdot)\colon\, I\to M^n$ are  integral curves of $T$ for any $x\in N^{m-1}$ and $\psi(\cdot, s)\colon\, N^{m-1}\to M^m$ are  leaves of $\{T\}^\perp$ for any $s\in I$.
Denoting by $E_1$ and $E_2$ the distributions given by tangent spaces to the leaves of the product foliation of  $N^{m-1}\times I$, we have that $E_1$ and $E_2$ are mutually orthogonal and $E_2$ is totally geodesic with respect to the metric induced by $\psi$. Set   $\tilde{f}=i\circ {f}\circ \psi$. Then 
\be \label{eq:tildef*}\<\tilde{f}_*X, i_*\frac{\d}{\d t}\>=\<\psi_*X, T\>=0\ee
for every $X\in E_1$. Moreover,  $\alpha_{\tilde{f}}(X,\frac{\d}{\d s})=0$ for every $X\in E_1$, in view of (\ref{eq:eigenvectorAxi}) and the fact that $f$ belongs to class ${\cal A}$. Hence, using that $E_2$ is totally geodesic we obtain that 
$$\tilde{\nabla}_{\frac{\d}{\d s}}\tilde{f}_*X=\tilde{f}_*\nabla_{\frac{\d}{\d s}}X+\alpha_{\tilde{f}}\left(X,\frac{\d}{\d s}\right)=\tilde{f}_*\nabla_{\frac{\d}{\d s}}X\in \tilde{f}_*E_1$$
for all $X\in E_1$, hence $\tilde{f}_*E_1$ is constant 
in $ \Ee^{n+2}$ along the leaves of $E_2$. In view of (\ref{eq:tildef*}) we can assume that $\tilde{g}:=\tilde{f}(\cdot , 0)$ satisfies $\tilde{g}(N^{m-1})\subset \Q_\epsilon^n\times \{0\}$. Set $\hat f=\pi\circ \tilde f$ and $h=\<\tilde f, i_*\d/\d t\>$, so that 
$$\tilde f=\hat f+hi_*\frac{\d}{\d t}.$$
Using (\ref{eq:tildef*}) we obtain
$$\<\tilde{f}_*X, \tilde{f}\>=\<\hat f_*X, \hat f\>=0$$
for all $X\in E_1$, since $\<\hat f, \hat f\>=\epsilon$. Therefore, we have
$$\tilde f(x, s)\in (\tilde{f}_*(x,s)E_1(x, s))^\perp=(\tilde{f}_*(x,0)E_1(x, 0))^\perp=(\tilde{g}_*(x)T_xN)^\perp,$$
where in the first equality we have used that  $\tilde{f}_*E_1$ is constant 
in $\Ee^{n+2}$ along $E_2$. Hence, for fixed $s\in I$, we have that $\xi_s(x):=\tilde{f}(x,s)$ defines a normal vector field along $\tilde{g}$.  Moreover,
$$\tilde{\nabla}_X \xi_s\in \tilde{f}_*(x,s)E_1(x,s)=\tilde{g}_*(x)T_xN,$$
thus $\xi_s$ is parallel along $\tilde{g}$ in the normal connection. It follows that $$x\in N\mapsto\spa\{\tilde{f}(x,s)\,:\,s\in I\}$$ is a parallel flat subbundle of $N_{\tilde g}N$ and, for fixed $x\in N$, the fiber $\{\tilde{f}(x,s)\,:\,s\in I\}$
is contained in a cylinder $\Q_\epsilon \times I\subset N^{\tilde{g}}_{x}N$, for $\<\hat{f}(x,s), \hat{f}(x, s)\>=\epsilon$.

     Let $g\colon\,N^{m-1}\to \Q_\epsilon^n$  be defined by $\tilde{g}=k\circ g$, and let $\{\xi_1, \ldots, \xi_k\}$ be an orthonormal set of parallel normal vector fields along $g$ such that $\tilde{\xi}_i=i_*\xi_i$, $1\leq i\leq k$,  $\tilde{\xi}_0=\tilde g$ and $ \tilde{\xi}_{k+1}=i_*\d/\d t$ span 
$\{\tilde{f}(x,s)\,:\,s\in I\}$ for each $x\in N^{m-1}$. Note that 
$$X\<\tilde{f}, \tilde{\xi}_i\>=\<\tilde{f}_*X,\tilde{\xi}_i\>+\<\tilde{f},\tilde{\nabla}_X\tilde{\xi}_i\>=0,$$
for $\tilde{f}$ is a normal vector field, $\tilde{f}_*(x,s)X\in \tilde{g}_*(x)T_xN$ and $\tilde{\xi}_i$ is parallel in the normal connection of $\tilde{g}$. Then we can write $$\tilde{f}(x,s)=\sum_{i=0}^{k+1}\alpha_i\tilde{\xi}_i,\,\,\,\mbox{with}\,\,\,\alpha_i=\alpha_i(s).$$
Moreover, from $\<\hat f, \hat f\>=\epsilon$ we obtain that $\epsilon \alpha_0^2+\sum_{i=1}^k\alpha_i^2=\epsilon$.\qed
\vspace{2ex}

\noindent \emph{Proof of Corollary \ref{cor:fnb}:}  It follows from the Ricci equation (\ref{eq:ricci2}) that $f$ has flat normal bundle  if and only if all shape operators $A^f_\zeta$, $\zeta\in N^fM$, are simultaneously diagonalizable, whereas $\tilde f$ has flat normal bundle if and only if this holds for all shape operators $A^{\tilde f}_\xi$, $\xi\in N^{\tilde{f}} M$. Since the vector field $T$ in (\ref{eq:ddt}) is nowhere vanishing, the equivalence between $(i)$ and $(ii)$ then follows from  Proposition \ref{prop:xi} and the fact that $A^{\tilde f}_{i_*\zeta}=A^f_\zeta$ for all $\zeta \in N^fM$. 

 Let $f\colon\, M^m\to \Q_\e^n\times \R$ be an isometric immersion in class ${\cal A}$. By Theorem \ref{thm:main}, it is locally given as in (\ref{eq:f}) in terms of an isometric immersion $g\colon\, N^{m-1}\to \Q_\e^n$. Since $T$ is an eigenvector of all shape operators of $f$, it follows from  (\ref{eq:alphapt}) and (\ref{eq:sffs}) that all shape operators of $f$ commute if and only if the same holds for the shape operators of $g$. By the Ricci equation, we conclude that $ f$ has flat normal bundle if and only if the same holds for $g$. Hence $(ii)$ and $(iii)$ are equivalent.  \qed\vspace{2ex}
 
\noindent \emph{Proof of Corollary \ref{cor:constantangle}:} Let $f\colon\, M^m:=N^{m-1}\times I\to \Q_\e^n\times \R$ be given  by
(\ref{eq:f}) with $\alpha\colon\, I\to \Q_\e^k\times \R$  a geodesic of $\Q_\e^k\times \R$. Then $f$ belongs to class ${\cal A}$ by Theorem \ref{thm:main} and $\alpha_f(T,T)=0$ by (\ref{eq:alphapt3}).  Thus $\eta$ is parallel in the normal connection of $f$ by (\ref{eq:Derivadaeta}).

Conversely, let $f\colon\, M^m:=N^{m-1}\times I\to \Q_\e^n\times \R$ be an isometric immersion with the property that the vector field  $\eta$ is parallel in the normal connection.   We obtain  from (\ref{eq:Derivadaeta}) that $f$ belongs to   class ${\cal A}$  and that $\alpha_f(T,T)=0$. By Theorem \ref{thm:main}, it is locally given  by
(\ref{eq:f}) in terms of an isometric immersion $g\colon\,N^{m-1}\to \Q_\epsilon^n$ and a smooth regular curve  $\alpha\colon\, I\to \Q_\e^k\times \R$. That $\alpha$ is a geodesic follows from $\alpha_f(T,T)=0$ and   (\ref{eq:alphapt3}). 
\qed

\section{Rotational submanifolds in $\Q_\e^n\times\R$}

In this section we define rotational submanifolds in $\Q_\e^n\times\R$ with curves as profiles, extending the definition in \cite{dfv} for the hypersurface case. Then we prove Corollary \ref{prop:rot} in the introduction.

Let $(x_0,\ldots,x_{n+1})$ be the standard  coordinates on $\Ee^{n+2}$  with respect to which the flat metric  is written as 
$$d s^2=\e \,d x_0^2+d
x_2^2+\ldots+d x_{n+1}^2.$$
Regard $\Ee^{n+1}$ as 
$$\Ee^{n+1}=\{(x_0,\ldots,x_{n+1})\in\Ee^{n+2}:x_{n+1}=0\}$$
and 
$$ \Q_\e^n=\{(x_0,\ldots,x_{n})\in\Ee^{n+1}:\e\,x_0^2+x_2^2+\ldots+x_{n}^2=\e\}\,\,\,\mbox{( $x_0>0$ if $\e=-1$)}.$$

Let $P^{n-m+3}$ be a  subspace of $\Ee^{n+2}$ of dimension $n-m+3$ containing the $e_0$ and the 
$e_{n+1}$ directions, where $\{e_0, \ldots, e_{n+1}\}$ is the canonical basis. Then $$(\Q_\e^n\times\R)\cap P^{n-m+3}=
\Q_\e^{n-m+1}\times\R.$$  Denote by $\mathcal I$
the group of isometries of $\Ee^{n+2}$ that fix pointwise a  subspace $P^{n-m+2}\subset P^{n-m+3}$  also
containing the $e_{n+1}$ direction. 
Consider a curve  $\alpha$  in
$\Q_\e^{n-m+1}\times\R\subset P^{n-m+3}$ that lies in one of the two half-spaces of $P^{n-m+3}$ determined by  $P^{n-m+2}$.

\begin{definition} {\em An $m$-dimensional {\em rotational submanifold  in $\Q_\e^n\times\R$ with profile curve $\alpha$ and axis $P^{n-m+2}$}
is the orbit of $\alpha$ under the action of $\mathcal I$.}
\end{definition}

We will always assume that $P^{n-m+3}$ is spanned by $e_0, e_m,\ldots, e_{n+1}$.
In the case $\e=1$, we also assume that $P^{n-m+2}$ is spanned by $e_m,\ldots, e_{n+1}$.
Writing the curve $\alpha$ as 
$$\alpha(s)=\alpha_0(s)e_0+\sum_{i=m}^{n}\alpha_{i-m+1}(s)e_i+h(s)e_{n+1},$$
with $\sum_{i=0}^{n-m+1}\alpha_i^2=1$,
 the  rotational submanifold   in $\Sf^n\times\R$ with profile curve $\alpha$ and axis $P^{n-m+2}$ can be parametrized by
\be\label{eq:parametrization} \tilde{f}(s,t)=(\alpha_0(s)\varphi_1(t),\ldots,\alpha_0(s)\varphi_{m}(t),\alpha_1(s), \ldots, \alpha_{n-m+1}(s),h(s)),\ee
where $t=(t_1,\ldots,t_{m-1})$ and 
$\varphi=(\varphi_1,\ldots,\varphi_m)$  parametrizes
 $\Sf^{m-1}\subset\R^m$. 
 
 For $\e=-1$, one has three distinct possibilities, according as $P^{n-m+2}$ is Lorentzian, Riemannian or degenerate, respectively, and the rotational submanifold is called accordingly of {\em spherical\/}, {\em hyperbolic\/} or {\em parabolic\/} type. 
 In the first case, we can assume that 
 $P^{n-m+2}$ is spanned by $e_0, e_{m+1},\ldots, e_{n+1}$
and that 
\be\label{eq:paramcurve}\alpha(s)=\alpha_0(s)e_0+\sum_{i=m}^{n}\alpha_{i-m+1}(s)e_i+h(s)e_{n+1},\ee
with $-\alpha_0^2(s)+\sum_{i=1}^{n-m+1}\alpha_i^2=-1$. 
Then, the submanifold can be parametrized by
$$ \tilde{f}(s,t)=(\alpha_0(s),\alpha_{1}(s)\varphi_1(t),\ldots,\alpha_{1}(s)\varphi_{m}(t),\alpha_{2}(s),\ldots, \alpha_{n-m+1}(s), h(s)),$$
where again $t=(t_1,\ldots,t_{m-1})$ and 
$\varphi=(\varphi_1,\ldots,\varphi_m)$  parametrizes
 $\Sf^{m-1}\subset\R^m$.

In the second case, we can assume that 
 $P^{n-m+2}$ is spanned by $e_m, \ldots, e_{n+1}$.
 Then, with   the curve $\alpha$ also given as in  (\ref{eq:paramcurve}),  a  parametrization is
$$ \tilde{f}(s,t)=(\alpha_0(s)\varphi_1(t),\ldots,\alpha_0(s)\varphi_m(t),\alpha_1(s),\ldots, \alpha_{n-m+1}(s), h(s)),$$
where $t=(t_1,\ldots,t_{m-1})$ and
$\varphi=(\varphi_1,\ldots,\varphi_m)$  parametrizes
 $\Hy^{m-1}\subset\Le^m$.

 Finally, when $P^{n-m+2}$ is degenerate, we choose a pseudo-orthonormal basis 
 $$\hat{e}_0=\frac{1}{\sqrt{2}}(-e_0+e_{n}),\,\,\,\, \hat{e}_{n}=\frac{1}{\sqrt{2}}(e_{0}+e_{n}), \,\,\,\hat{e}_j=e_{j}, $$
 for $j\in \{1,\ldots, n-1,n+1\}$, and assume that   $P^{n-m+2}$ is spanned by $\hat{e}_{m}, \ldots, \hat{e}_{n+1}$.
 %, \hat{e}_{n}$ and $\hat{e}_{n+1}$. 
 Notice that $\<\hat{e}_0,\hat{e}_0\>=0=\<\hat{e}_{n},\hat{e}_{n}\>$ and $\<\hat{e}_0,\hat{e}_{n}\>=1$. Then, we can parametrize $\alpha$ by 
 $$\alpha(s)=\alpha_0(s)\hat{e}_0+\sum_{i=m}^{n}\alpha_{i-m+1}(s)\hat{e}_i+h(s)\hat{e}_{n+1},$$
with $2\alpha_0(s)\alpha_{n-m+1}(s)+\sum_{i=1}^{n-m}\alpha_i^2(s)=-1$,
 and a parametrization of the rotational submanifold  is 
  \be\label{eq:parametrization4} \tilde{f}(s, t)=(\alpha_0,\alpha_0t_1,\ldots, \alpha_0t_{m-1},\alpha_1, \ldots, \alpha_{n-m}, \alpha_{n-m+1}-\frac{\alpha_0}{2}\sum_{i=1}^{m-1}t_i^2, h),\ee
 where $t=(t_1, \ldots, t_{m-1})$ parametrizes $\R^{m-1}$, $\alpha_i=\alpha_i(s)$, $0\leq i\leq n-m+1$, and $h=h(s)$.\vspace{2ex}
 
 \noindent \emph{Proof of Corollary \ref{prop:rot}:}
 We can write (\ref{eq:parametrization}) as 
 $$\tilde{f}(s,t)=\alpha_0(s)\hat g(t)+\sum_{i=m}^{n}\alpha_{i-m+1}(s)e_i+h(s)e_{n+1},$$
 where $\hat g(t)=\sum_{i=1}^m\varphi_i(t)e_i$ for $t=(t_1, \ldots, t_{m-1})$. This shows that for  $\e=1$  a rotational submanifold is given as in (\ref{eq:f})  in terms of a totally geodesic  isometric immersion $\hat g\colon\,\Sf^{m-1}\to \Sf^n$.
 %, with the vector subbundle $E$ being the normal bundle $N^{\hat g}\Sf^{m-1}$. 
 The case of a rotational submanifold of hyperbolic type in $\Hy^n\times \R$ is similar. In particular, this proves that $(i)$ implies $(ii)$ in these cases. 
 
    Equation (\ref{eq:parametrization4}) can be written  as 
 \be\label{eq:p5} \tilde{f}(s, t)=\alpha_0\hat{g}+\sum_{i=m}^{n}\alpha_{i-m+1}(s)\hat{e}_i+h(s)\hat{e}_{n+1},\ee
    where $\hat{g}(t)=\hat{e}_0+\sum_{i=1}^{m-1}t_i\hat{e}_i-\frac{1}{2}(\sum_{i=1}^{m-1}t_i^2)\hat{e}_n.$ Notice that $\hat{g}$ defines an isometric immersion of $\R^{m-1}$ into $\Le^{n+2}$ (in fact into the light-cone $\Ve^{n+1}$, for $\<\hat{g}, \hat{g}\>=0$), and that $\hat{g},\hat{e}_m, \ldots, \hat{e}_n, \hat{e}_{n+1}$ is a pseudo-othonormal basis of $N^{\hat{g}}\R^{m-1}$, with  $\<\hat{g}, \hat{g}\>=0=\<\hat{e}_n, \hat{e}_n\>$, $\<\hat{g},e_n\>=1$ and $\hat{e}_m, \ldots, \hat{e}_{n-1}, \hat{e}_{n+1}$ an orthonormal basis of $\spa\{\hat{g}, e_n\}^\perp$. For any  fixed $s_0\in I$, let ${g}\colon\, \R^{m-1}\to \Hy^{n}$ be given by
     where  $v=\sum_{i=m}^{n}\alpha_{i-m+1}(s_0)\hat{e}_i$.
     Then ${g}$ defines an umbilical isometric immersion with the same normal space in $\Le^{n+2}$ as $\hat{g}$ at every $t\in \R^{m-1}$, i.e., 
    $$\spa\{\hat g, \hat{e}_m, \ldots, \hat{e}_n, \hat{e}_{n+1}\}=\spa\{g,\tilde{\xi}_1, \ldots, \tilde{\xi}_{n-m+1}, \hat{e}_{n+1}\},$$
    where $\tilde{\xi}_i=i_*\xi_i$, $1\leq i\leq n-m+1$, for a parallel orthonormal frame $\xi_1, \ldots, \xi_{n-m+1}$ of $N^{{g}}\R^{m-1}$.  Hence we can also write (\ref{eq:p5}) as
    $$\tilde{f}(s, t)=\tilde{\alpha}_0{g}+\sum_{i=1}^{n-m+1}\tilde{\alpha}_{i}(s)\tilde{\xi}_i+h(s)\hat{e}_{n+1},$$
     where $\tilde \alpha\colon\, I\to \Ee^{n-m+3}$ is a regular curve satisfying $-\tilde{\alpha}_0^2+\sum_{i=1}^{n-m+1}\tilde{\alpha}_i^2=-1$. Thus condition $(ii)$ holds for $f$. The case of a spherical rotational submanifold is similar and easier.
     
    Now suppose that $f$ is given  as in (\ref{eq:f})  in terms of an umbilical isometric immersion $g\colon\,N^{m-1}\to \Q_\epsilon^n$ (a geodesic circle if $m=2$).
    %, with the vector subbundle $E$ being the normal bundle $N^gN$. 
    Suppose first that $\e=1$. We can assume that the affine hull of  $g(N^{m-1})$ in $\R^{n+1}$ is $v+W$, where  $W$ is  the subspace spanned by $\{e_0, \ldots, e_{m-1}\}$ and  $v\in W^\perp$, hence $g=a\hat g+ v$, where $a\in \R$  and $\hat g$ is the composition $\hat g=i\circ \tilde g$ of a homothety $\tilde g\colon\, N^{m-1}\to \Sf^{m-1}$ with the canonical inclusion $i$ of $\Sf^{m-1}$ into $W=\R^m$ as the unit sphere centered at the origin. Then $g$ and $\hat g$ have the same normal spaces in $\R^{n+1}$ at every point of $N^{m-1}$, that is,  
    $$\spa\{\hat g, \hat{e}_m, \ldots, \hat{e}_n, \hat{e}_{n+1}\}=\spa\{g,\tilde{\xi}_1, \ldots, \tilde{\xi}_{n-m+1}, \hat{e}_{n+1}\},$$
    where $\tilde{\xi}_i=i_*\xi_i$, $1\leq i\leq n-m+1$, for a parallel orthonormal frame $\xi_1, \ldots, \xi_{n-m+1}$ of $N^{{g}}\Sf^{m-1}$.  Hence $f$ can also be parametrized by 
        $$\tilde{f}(s, t)=\tilde{\alpha}_0\hat{g}+\sum_{i=m}^{n}\tilde{\alpha}_{i-m+1}(s)\hat{e}_i+h(s)\hat{e}_{n+1},$$
     where $\tilde \alpha\colon\, I\to \Ee^{n-m+3}$ is a smooth regular curve satisfying $\sum_{i=0}^{n-m+1}\tilde{\alpha}_i^2=1$. Thus $f$ is a rotational submanifold with $\tilde{\alpha}$ as profile.
     
     If $\e=-1$, we argue for the parabolic case, the others being similar and easier. We can assume that $N^{m-1}=\R^{m-1}$ and that the the affine hull of $g(\R^{m-1})$ in $\Le^{n+1}$ is $v+W$, where  $W$ is the subspace spanned by $\{\hat{e}_1, \ldots, \hat{e}_{m-1}, \hat{e}_{n}\}$ and $v\in W^\perp$. Then $g=a\hat g+ v$, where $a\in \R$ and $\hat{g}(t)=\hat{e}_0+\sum_{i=1}^{m-1}t_i\hat{e}_i-\frac{1}{2}(\sum_{i=1}^{m-1}t_i^2)\hat{e}_n$ for 
 $t=(t_1, \ldots, t_{m-1})$. As before, by using the fact that $g$ and $\hat g$ have the same normal spaces in $\Le^{n+1}$ for every $t\in \R^{m-1}$, we conclude that $f$ is a rotational submanifold parametrized as in (\ref{eq:p5}).

 The second fundamental form of $f$ being given by (\ref{eq:sff0}) is equivalent to the restriction of each shape operator $A_\xi^f$ to $\{T\}^\perp$ being a multiple of the identity tensor. In particular, if it is satisfied then the immersion $f$ is in class ${\cal A}$, hence it is locally given as in (\ref{eq:f}) in terms of an  isometric immersion $g\colon\, N^{m-1}\to \Q_\e^n$. It follows from formulas (\ref{eq:alphapt})-(\ref{eq:sffs}) in part $(iv)$ of Proposition \ref{prop:cu} that $g$ is umbilical.  Conversely, if $f$ is locally given as in (\ref{eq:f}) in terms of an umbilical isometric immersion $g\colon\, N^{m-1}\to \Q_\e^n$, then formulas (\ref{eq:alphapt})-(\ref{eq:sffs}) imply that the restriction of each shape operator $A_\xi^f$ to $\{T\}^\perp$ is a multiple of the identity tensor, hence the second fundamental form of $f$ is as in (\ref{eq:sff0}). To conclude the proof that $(ii)$ and $(iii)$ are equivalent, it remains to show that for  $m=2$ the additional assumption that the vector field $\zeta$ in (\ref{eq:sff0}) be parallel along $\{T\}^\perp$ is equivalent to the unit-speed curve  $g\colon\,J:=N^{m-1}\to \Q_\e^n$ being a  geodesic circle. 
 
 Write $\zeta=\alpha_f(X,X)$, where $X$ is a unit vector field orthogonal to $T$. Let  $\tilde{f}=i\circ f$. In view of (\ref{eq:eigenvectorAxi}) we have
 $$\alpha_{\tilde f}(X,X)=i_*\alpha_f(X,X)-\nu,$$
 hence we obtain using (\ref{eq:nconns}) and (\ref{eq:normalnu}) that
 \begin{eqnarray*}
 \tilde \nabla_X \alpha_{\tilde f}(X,X)&=&-\tilde{f}_*A^{\tilde f}_{\alpha_{\tilde f}(X,X)}X+\tilde\nabla_X^\perp\alpha_{\tilde f}(X,X)\\
 &=&-\tilde{f}_*A^{\tilde f}_{\alpha_{\tilde f}(X,X)}X+i_*\nabla_X^\perp\alpha_f(X,X).
 \end{eqnarray*}
 Thus, $\nabla^\perp_X\alpha_f(X,X)=0$ if and only if 
 \be\label{eq:condsurf}\tilde \nabla_X \alpha_{\tilde f}(X,X)=-\tilde{f}_*A^{\tilde f}_{\alpha_{\tilde f}(X,X)}X.\ee
 It follows from  (\ref{eq:alphapt}) that at the point $(t,s)$ we have 
 \be\label{eq:alphatil} \alpha_{\tilde f}(X,X)=\frac{\tilde{g}''(t)}{\<\tilde{g}''(t),\phi_t(\bar{\alpha}(s))\>},\ee
 where  $\tilde g=k\circ g$. 
 From (\ref{eq:diff1}) we obtain
 \be\label{eq:X}X=\frac{1}{\<\tilde{g}''(t),\phi_t(\bar{\alpha}(s))\>}\frac{d}{dt}^{\cal H},\ee
 where $\frac{d}{dt}$ is a unit vector field along $J$. Hence
 $$\tilde \nabla_X \alpha_{\tilde f}(X,X)=-\frac{\<\tilde{g}'''(t),\phi_t(\bar{\alpha}(s)\>}{\<\tilde{g}''(t),\phi_t(\bar{\alpha}(s))\>^3}\tilde{g}''(t)+
 \frac{1}{\<\tilde{g}''(t),\phi_t(\bar{\alpha}(s))\>^2}\tilde{g}'''(t).$$
 On the other hand, equations (\ref{eq:diff1}), (\ref{eq:alphatil}) and (\ref{eq:X})  yield
 $$\tilde{f}_*A^{\tilde f}_{\alpha_{\tilde f}(X,X)}X=\frac{\<\tilde{g}''(t),\tilde{g}''(t)\>}{\<\tilde{g}''(t),\phi_t(\bar{\alpha}(s))\>^2}\tilde{g}'(t).$$
 It follows easily that (\ref{eq:condsurf}) holds if and only if 
 $$\tilde{g}'''(t)=-\<\tilde{g}''(t),\tilde{g}''(t)\>\tilde{g}'(t),$$
 which is equivalent to $g$ being  a geodesic circle. \qed

\section{Reduction of codimension}

In this section we prove Lemma \ref{le:redcod} and Theorem \ref{thm:dajczer} stated in the introduction.\vspace{2ex}\\
\noindent\emph{Proof of Lemma \ref{le:redcod}:} We have  from (\ref{eq:Derivadaeta}) that $\nabla_X^\perp\eta\in N_1\subset L$ for every $X\in TM$. Since $\nabla^\perp N_1\subset L$ by assumption, it follows that $L$ is a parallel subbundle of $N^fM$.  Let  $\tilde{f}=i\circ f$, where $i\colon\,\Q_\e^n\times\R\to  \Ee^{n+2}$ is the  inclusion.  
Given $\xi\in L^\perp=N_1^\perp\cap \{\eta\}^\perp$, 
  from (\ref{eq:nconns}) and the fact that $L$ is a parallel subbundle of $N^fM$ we obtain
 $$\tilde \n^\perp_Xi_*\xi=i_*\n^\perp_X\xi\in i_*L^\perp,$$
 hence $i_*L^\perp$ is a parallel subbundle of $N^{\tilde f}M$.

 Since $i_*L^\perp\subset\tilde {N}_1^\perp$, where $\tilde N_1(x)$ is the first normal space of $\tilde f$ at $x\in M^m$, it follows that $i_*L^\perp$ is a constant subspace of $\Ee^{n+2}$, which is orthogonal to  $\frac{\d}{\d t}$. Denote by $K$ the orthogonal complement of $i_*L^\perp$ in $\Ee^{n+2}$. 
 Then, for any fixed $x_0\in M^m$ we have  
 $$\tilde{f}(M^m)\subset \tilde{f}(x_0)+K.$$
 But since $K$ contains  $\frac{\d}{\d t}$ and $\nu(x_0)$, it also contains the position vector $\tilde{f}(x_0)$. Thus $\tilde{f}(x_0)+K=K$. 
 We conclude that 
$\tilde f(M)\subset (\Q_\e^n\times \R)\cap K=\Q_\e^{m+\ell-1}\times \R.$\vspace{2ex}\qed\\
 \noindent\emph{Proof of Theorem \ref{thm:dajczer}:}
 Assume that $\nabla^\perp N_1\subset L$. Then condition $(ii)$ is trivially satisfied. To prove $(i)$, 
 first notice that for $\xi \in N_1^\perp$ the Ricci equation gives
 $$R^\perp(X,Y)\xi=\alpha(X,A_\xi Y)-\alpha(A_\xi X, Y)=0.$$
 Given  $\xi\in L^\perp$, we have that  $\xi\in N_1^\perp$ and that $\nabla^\perp_Z\xi\in N_1^\perp$ by our assumption, hence
$$(\nabla_ZR^\perp)(X,Y,\xi)=\nabla_ZR^\perp(X,Y)\xi-R^\perp(\nabla_ZX,Y)\xi-R^\perp(X,\nabla_ZY)\xi-R^\perp(X,Y)\nabla_Z^\perp\xi=0.$$

   To prove the converse,  let $\xi\in L^\perp$. Since $R^\perp(X,Y)\xi=0$ for all $X,Y\in TM$, we obtain from $(i)$ that 
   $$R^\perp(X,Y)\nabla_Z^\perp\xi=0$$
   for all $X,Y,Z\in TM$. Using the Ricci equation again, we obtain that
   $$[A_{\nabla_Z^\perp \xi}, A_{\nabla_W^\perp \xi}]=0$$
   for all $Z, W\in TM$. Hence, at any $x\in M$ there exists an orthonormal basis $Z_1, \ldots, Z_n$ of $T_xM$ that 
   diagonalizes simultaneously all shape operators $A_{\nabla_Z^\perp \xi}$, $Z\in TM$.   We will show that
   $$\<\nabla^\perp_{Z_k}\xi, \alpha(Z_i, Z_j)\>=0$$
    for all $1\leq i,j,k\leq n$, which implies that $\nabla^\perp_X\xi\in N_1^\perp$ for all $X\in TM$.  
    
    From the choice of the basis $Z_1,\ldots, Z_n$, we have
    $$\<\alpha(Z_i, Z_j), \nabla_{Z_k}^\perp \xi\>=\<A_{\nabla^\perp_{Z_k}\xi}Z_i, Z_j\>=0$$
    if $i\neq j$.  It follows from the Codazzi equation (\ref{eq:codazzi3}) 
    and the fact that $\xi\in L^\perp\subset \{\eta\}^\perp$ that
    $$A_{\nabla_{Z_i}^\perp\xi}Z_k=A_{\nabla_{Z_k}^\perp\xi}Z_i,$$
    hence the eigenvalue $\lambda_{ki}$ of $A_{\nabla_{Z_k}^\perp\xi}$ correspondent to $X_i$ vanishes unless $k=i$. Therefore, 
    $$\<\alpha(Z_i, Z_i), \nabla_{Z_k}^\perp \xi\>=\<A_{\nabla_{Z_k}^\perp\xi}Z_i, Z_i\>=0, \,\,\,\mbox{if}\,\,i\neq k.$$
    Finally, the assumption $\nabla^\perp H\in L$ and the above imply that
    $$\<\alpha(Z_i, Z_i), \nabla_{Z_i}^\perp \xi\>=n\<H,\nabla_{Z_i}^\perp \xi\>=0.\qed$$
    
 \subsection{Alencar--do Carmo--Tribuzzi Theorem}
   
   In this subsection we apply Theorem \ref{thm:dajczer} to give a simple proof of the following  theorem due to  Alencar, do Carmo and Tribuzzi \cite{adct}.
   
   \begin{theorem}\po\label{thm:adct} Let $f\colon\,M^2\to \Q_\epsilon^n\times \R$, $n\geq 5$, be a surface with nonzero parallel mean curvature 
   vector. Then, one of the following possibilities holds:
   \begin{itemize}
   \item[(i)] $f$ is a minimal surface of a  umbilical hypersurface of a slice $\Q_\epsilon^n\times \{t\}$.
   \item[(ii)] $f$ is a surface with constant mean curvature in a three-dimensional  umbilical or totally geodesic submanifold of
  a slice $\Q_\epsilon^n\times \{t\}$.  
   \item[(iii)]  $f(M^2)$ lies in a totally geodesic submanifold $\Q_\epsilon^m\times \R$, $m\leq 4$, of $\Q_\epsilon^n\times \R$.
   \end{itemize}
   \end{theorem}
   \proof Since the mean curvature vector $H$ is parallel and nonzero, the function   $\mu:=\|H\|^2$ on $M^2$ is a nonzero constant. Suppose first that $A_H=\mu I$ everywhere on $M^2$. We claim that the vector field $T$ vanishes identically. Assuming otherwise, there exists an open subset $U$ where $T\neq 0$. Choose a unit vector field $X$ on $U$ orthogonal to $T$. Then
   \be\label{eq:h1}\<H, \alpha(X,T)\>=\mu\<X,T\>=0.\ee
   By the Codazzi equation (\ref{eq:codazzi2}) we have 
   $$\<(\nabla^\perp_T \alpha)(X,X)-(\nabla^\perp_X \alpha)(T,X),H\>=-\|T\|^2\<\eta, H\>.$$
   It follows easily from (\ref{eq:h1}) and the fact that $\mu$ is constant on $M^2$ that the left-hand-side of the 
   preceding equation is zero. Thus $\<\eta,H\>$ vanishes on $U$, and hence
   $$0=T\<\eta,H\>=\<\nabla_T^\perp\eta, H\>=-\<\alpha(T,T),H\>=-\mu\|T\|^2,$$ 
   where we have used (\ref{eq:Derivadaeta}) in the third equality. This is a contradiction and proves the claim.
   
   Therefore, if $A_H=\mu I$ everywhere on $M^2$ then $f(M^2)$ is contained in a slice $\Q_\epsilon^n\times \{t\}$ of $\Q_\epsilon^n\times \R$, and either of possibilities $(i)$ or $(ii)$ holds by  Theorem $4$ in \cite{yau}.

 Assume now that  $A_H\neq \mu I$ on an open  subset $V$ of $M^2$. 
 Since $H$ is parallel, it follows from the Ricci equation that $[A_H, A_\zeta]=0$ for any $x\in M^2$ and every normal vector $\zeta \in N_xM$. Then, the fact that $A_H$ has distinct eigenvalues on $V$ implies that the eigenvectors of $A_H$ are also eigenvectors of $A_\zeta$ for any $\zeta\in N_xM$, $x\in V$. Hence all shape operators are simultaneously diagonalizable at any $x\in V$, which implies that $f$ has flat normal bundle on $V$ by the Ricci equation (\ref{eq:ricci2}). In particular, the first normal spaces $N_1$ of $f$ have dimension at most two at any $x\in V$. Let $W\subset V$ be an open subset where $L=\dim N_1 + \spa\{\eta\}$ has constant dimension $\ell\leq 3$. It follows from Theorem \ref{thm:dajczer} that $f(W)$ lies in a totally geodesic submanifold $\Q_\epsilon^{2+\ell-1}\times \R$ of $\Q_\epsilon^n\times \R$. By analyticity of $f$ (see Remark $1$ of \cite{adct}), we conclude that $f(M^2)\subset \Q_\epsilon^{2+\ell-1}\times \R$.\vspace{2ex}\qed

 \noindent\emph{Proof of Corollary \ref{cor:pmc}:} Let $X$ be a unit vector field orthogonal to $T$. By Corollary \ref{prop:rot}, in order to prove that $f$ is a rotational surface it 
suffices to show that   $\nabla^\perp_X\alpha_f(X,X)=0$. We follow essentially the proof of Proposition $2$ in \cite{adct}. Since the mean curvature vector field 
$$H=\frac{1}{2}(\alpha_f(X,X)+\|T\|^{-2}\alpha_f(T,T))$$
is parallel in the normal connection, we have
$$\nabla^\perp_X\alpha_f(X,X)=-\nabla^\perp_X(\|T\|^{-2}\alpha_f(T,T))=-X(\|T\|^{-2})\alpha_f(T,T)-
\|T\|^{-2}\nabla^\perp_X\alpha_f(T,T).$$
Now, since $f$ is in class ${\cal A}$, we have from (\ref{eq:ntt}) that
$$\<\nabla_TT,X\>=0=\<\nabla_XT,T\>.$$
In particular, $X(\|T\|^{-2})=0$.
Moreover, using the Codazzi equation (\ref{eq:codazzi2}) we obtain
$$\begin{array}{l}\nabla^\perp_X\alpha_f(T,T)=(\nabla^\perp_X\alpha_f)(T,T)+2\alpha(\nabla_XT,T)=(\nabla^\perp_X\alpha_f)(T,T)=(\nabla^\perp_T\alpha_f)(X,T)\vspace{2ex}\\\hspace*{12.8ex}=\nabla^\perp_T\alpha_f(X,T)-\alpha_f(\nabla_TX,T)-\alpha_f(X,\nabla_TT)=0.
\end{array}$$
  That $f(M^2)$ is contained in a totally geodesic submanifold $\Q_\epsilon^m\times \R$, $m\leq 4$, and hence that its profile curve lies in a totally geodesic submanifold $\Q_\epsilon^s\times \R$, $s\leq 3$, follows from Theorem~\ref{thm:adct}.\qed
   
   \section{Umbilical submanifolds of $\Sf^n\times \R$}
   
  We are now in a position to prove Theorem \ref{thm:main2} in the introduction.\vspace{1ex}\\ 
 \noindent \emph{Proof of Theorem \ref{thm:main2}:} Since $\phi\colon\,\Sf^{m+1}\times \R\to \R^{m+2}\setminus \{0\}$  given by $\phi(x,t)=e^tx$ is a conformal diffeomorphism, it follows that $M^m_{r,h}=\phi^{-1}(\Sf^m_{r,h})$ is an umbilical submanifold of $\Sf^{m+1}\times \R$, for a conformal diffeomorphism preserves umbilical submanifolds. Assertion $(i)$ and completeness of $M^m_{r,h}$ are  clear, for $M^m_{r,h}=\phi^{-1}(\Sf^m_{r,h})$ if $(r,h)\neq (d,0)$ and $M_{d,0}=\phi^{-1}(\Sf^m_{d,0}\setminus \{0\})$. It is easily seen that the totally geodesic hypersurfaces $\Sf^{m}\times \R$ of $\Sf^{m+1}\times \R$ are the images by $\phi^{-1}$ of the hyperplanes through the origin in $\R^{m+2}$. Since
 $\Sf^m_{r,h}$ lies in such a hyperplane if and only if $h=0$, the assertion in $(ii)$ follows.  Assertion $(iii)$ follows from the fact that  $\Sf^m_{r,0}$ is homologous to zero in $\R^{m+1}$ if $r<d$ and inhomologous to zero in $\R^{m+1}$ if $r>d$. 
 
 We now prove $(iv)$. Since orthogonal transformations of $\R^{m+2}$ correspond under the diffeomorphism $\phi$ to isometries of $\Sf^{m+1}\times \R$ fixing pointwise the factor $\R$, and homotheties of $\R^{m+2}$ correspond to translations along $\R$, we can assume that $\bar x=\frac{\sqrt{2}}{2}(0,\ldots, 0,1)\in \R^{m+1}$. Let 
 $${\cal I}=\{(p,q)\in \R\,:\, (p-1)^2\leq q< p^2\}.$$ For each $(p,q)\in {\cal I}$, set $J_{p,q}=(-\sqrt{p-\sqrt{q}},\sqrt{p-\sqrt{q}})$ and define $h_{p,q}\colon\, \bar{J}_{p,q}\to \R$ by
 $$h_{p,q}(s)=\sqrt{p-s^2+\sqrt{(p-s^2)^2-q}}\,.$$
Let $Y_{p,q}\colon\, \Sf^{m-1}\times \bar{J}_{p,q}\to \Sf^{m+1}\times \R$ and  $Z_{p,q}\colon\, \Sf^{m-1}\times \bar{J}_{p,q}\to \Sf^{m+1}\times \R$ be given by
$$Y_{p,q}(x, s)=\left(sx, \frac{\sqrt{2}}{2}\left(h_{p,q}(s)+\frac{1-p}{h_{p,q}(s)}\right), \frac{\sqrt{q-(p-1)^2}}{\sqrt{2}\,h_{p,q}(s)}, \log h_{p,q}(s)\right)$$
and, for $q\neq 0$,
 $$Z_{p,q}(x, s)=\left(sx, \frac{\sqrt{2}}{2}\left(\frac{1-p}{\sqrt{q}}h_{p,q}(s)+\frac{\sqrt{q}}{h_{p,q}(s)}\right), \frac{\sqrt{q-(p-1)^2}}{\sqrt{2q}}h_{p,q}(s), \log \frac{\sqrt{q}}{h_{p,q}(s)}\right).$$
 Notice that $Z_{p,q}=\Psi\circ Y_{p,q}$, where $\Psi\colon\, \Sf^{m+1}\times \R\to \Sf^{m+1}\times \R$ is the isometry defined by 
 $\Psi(x,s)=(Ax,-s+\log \sqrt{q})$, with $A\in O(m)$ given by
 $$A=\left(\begin{array}{l}I_{m-2}\,\,\,\,\, 0\vspace{1ex}\\
 \,\,0\,\,\,\,\,\,\, \,\,\,B\end{array}\right), \,\,\,\,\,\,\,B=\frac{1}{\sqrt{q}}\left(\begin{array}{l}\,\,\,\,\,\,1-p\,\,\,\,\,\,\,\,\,\,\,\, \,\,\,\,\,\,\,\,\,\,\,\sqrt{q-(1-p)^2}\vspace{1ex}\\
 \sqrt{q-(1-p)^2}\,\,\, \,\,\,-(1-p)\end{array}\right).$$
 Let $\psi\colon\,{\cal I}\to (0, \infty)\times [0, \infty)$ be the diffeomorphism given by
 $$\psi(p,q)=\frac{\sqrt{2}}{2}\left(\sqrt{p^2-q}, \sqrt{q-(p-1)^2}\right).$$
 Then $(iv)$ is a consequence of the following fact.
 \begin{lemma}\po \label{prop:paramet} For $(r,h)=\psi(p,q)$ we have
 $$M^m_{r,h}=\left\{\begin{array}{l} Y_{p,q}(\Sf^{m-1}\times \bar{J}_{p,q})\cup Z_{p,q}(\Sf^{m-1}\times \bar{J}_{p,q}), \,\,\,\mbox{if}\,\,\,(r,h)\neq (d,0),\vspace{1ex}\\
 Y_{1,0}(\Sf^{m-1}\times (-1,1)),\,\,\,\,\,\,\mbox{if}\,\,\,(r,h)=(d,0). \end{array}\right.
 $$
 \end{lemma}
 \proof We argue for $(r,h)\neq (d,0)$, the case $(r,h)=(d,0)$ ( i.e., $(p,q)=(1,0)$) being similar and easier. A straightforward computation shows that
 $$(\phi\circ Y_{p,q})(x,s)=\left(sh_{p,q}(s)x, \frac{\sqrt{2}}{2}({h}^2_{p,q}(s)-p),0\right)+(\bar x, h)
 %\frac{\sqrt{2}}{2}\left(0,1,\sqrt{q-(p-1)^2}\right)
 $$
and $$(\phi\circ Z_{p,q})(x,s)=\left(s\bar{h}_{p,q}(s)x, \frac{\sqrt{2}}{2}(\bar{h}^2_{p,q}(s)-p),0\right)+
%\frac{\sqrt{2}}{2}\left(0,1,\sqrt{q-(p-1)^2}\right).
(\bar x, h),$$
where $\bar{h}_{p,q}(s)=\sqrt{q}/h_{p,q}(s)$. Let $\gamma\colon\, \bar{J}_{p,q}\to \R^2$ and $\bar \gamma\colon\, \bar{J}_{p,q}\to \R^2$ be given by $$\gamma(s)=(sh_{p,q}(s), \frac{\sqrt{2}}{2}(h^2_{p,q}(s)-p))\,\,\,\mbox{and}\,\,\,\bar{\gamma}(s)=(s\bar{h}_{p,q}(s), \frac{\sqrt{2}}{2}(\bar{h}^2_{p,q}(s)-p)),$$ respectively. 
 Then, the statement follows from the fact that $\gamma(\bar{J}_{p,q})\cup \bar\gamma(\bar{J}_{p,q})$  
is the circle of radius $r=\sqrt{\frac{p^2-q}{2}}$ centered at the origin.
\qed\vspace{1ex}
 
 We now prove the converse. Let $f\colon\, M^m\to \Sf^n\times \R$, $m\geq 2$, be an umbilical isometric immersion. If the vector field $T$ in (\ref{eq:ddt}) vanishes identically, then $f(M^n)$ is contained in a slice $\Sf^n\times \R$, and this gives the first possibility in the statement. Now assume that $T$ does not vanish at some point, and hence on some open subset $U\subset M^n$. It suffices to prove that there exist  open subsets $\tilde U\subset U$ and $V\subset \Sf^m$,  $(p,q)\in {\cal I}$ and an interval $I\subset I_{p,q}$  such that, up to an isometry of $\Sf^n\times \R$, we have  $f(\tilde U)\subset Y_{p,q}(V\times I)$. For  this implies that $f(\tilde U)\subset M^m_{r,h}$ with $(r,h)=\psi(p,q)$, and thus $(\phi\circ f)(\tilde U)\subset \phi(M^m_{r,h})=\Sf^m_{r,h}$. Since $\phi\circ f$ is an umbilical immersion into $\R^{n+2}\setminus \{0\}$, it follows that $(\phi\circ f)(M^m)\subset \Sf^m_{r,h}$, and hence $ f(M^m)\subset M^m_{r,h}$.
 
 From Codazzi equation (\ref{eq:codazzi2}) and $\alpha_f(X,Y)=n\<X, Y\>H$ for all $X, Y\in TM$, where $H$ is the mean curvature vector of $f$, we obtain
 \be\label{eq:nH}n\nabla_X^\perp H=-\epsilon\<X,T\>\eta\ee
 for every $X\in TM$. If $H$ and $\eta$ are linearly dependent on $U$, it follows from Lemma~\ref{le:redcod} that   $f$ has substantial codimension one on $U$. Otherwise, there exists an open subset $\tilde U\subset U$ such that
 $H$ and $\eta$ are linearly independent on $\tilde U$, in which case  Lemma \ref{le:redcod} implies that $f$ has substantial codimension two on $\tilde U$. 
 
 On the other hand, since $f$ is umbilical and its mean curvature vector $H$ is parallel in the normal connection along $\{T\}^\perp$ by (\ref{eq:nH}),  condition $(iii)$ in  Corollary \ref{prop:rot} is satisfied. Thus $f$ is a rotational submanifold. 
 
    Summing up,    $f|_{\tilde U}$ is a rotational submanifold
of substantial codimension at most two over a curve in a totally geodesic submanifold $\Sf^s\times \R$, $s\leq 2$. Hence, we can assume that $n=m+1$ and $s=2$. Equivalently, in view of the last assertion in Corollary \ref{prop:rot},  we obtain that $f|_{\tilde U}$ is 
given by (\ref{eq:f}) in terms of a totally geodesic isometric immersion $g\colon\,V\subset \Sf^{m-1}\to \Sf^{m}$ and a regular curve 
$\alpha\colon\, I\to \Sf^2\times \R\subset \R^4$, $\alpha=(\alpha_0, \alpha_1, \alpha_2, \alpha_3)$,  with $\alpha_0^2+\alpha_1^2+\alpha_2^2=1$. 

With notations as in Proposition \ref{prop:cu}, we  have by (\ref{eq:N}) and the umbilicity of $f$ that  $A_{\phi_x(\zeta)}^{\tilde f}$ is a multiple of the identity tensor for every $\zeta\in \alpha'(s)^\perp \cap \bar{\alpha}(s)^\perp$. Using that $P_s(x)=\alpha_0(s)I$, it follows from  (\ref{eq:alphapt}), (\ref{eq:alphapt3}) and (\ref{eq:sffs}) that 
$$-\frac{\<\tilde{g},\phi_x(\zeta)\>}{\alpha_0(s)}=\frac{\<\alpha''(s), \zeta\>}{\<\alpha'(s),\alpha'(s)\>}\,\,\,\,\mbox{for all}\,\,\,\,\zeta\in \alpha'(s)^\perp \cap \bar{\alpha}(s)^\perp,$$
or equivalently,
$$\<\alpha_0(s)\alpha''(s)+\varphi(s) e_0, \zeta\>=0\,\,\,\,\mbox{for all}\,\,\,\,\zeta\in \alpha'(s)^\perp \cap \bar{\alpha}(s)^\perp,$$
since $\tilde g=\phi_x(e_0)$. Here $\varphi(s)=\<\alpha'(s), \alpha'(s)\>$. Hence, there exist smooth functions $y=y(s)$ and $z=z(s)$ such that
\be\label{eq:umbcond}\alpha_0\alpha''+\varphi e_0=y\alpha'+z\bar{\alpha}.\ee
We write the preceding equation as 
\be\label{eq:umbcond2}\alpha_0\bar{\alpha}''+\alpha_0\alpha_3''e_3+\varphi e_0=y\bar{\alpha}'+y\alpha_3'e_3+z\bar{\alpha}.\ee
Notice that 
\be\label{eq:baralpha}
 \<\bar{\alpha}, \bar{\alpha}\>=1, \,\, \<\bar{\alpha}', \bar{\alpha}\>=0, \,\, \<\bar{\alpha}', \bar{\alpha}'\> = \varphi-(\alpha'_3)^2=-\<\bar{\alpha}'', \bar{\alpha}\>, \,\, \text{and} \,\, \<\bar{\alpha}'',\bar{\alpha}'\>=\frac{1}{2}(\varphi'-2\alpha_3'\alpha_3'').
\ee
On the other hand, taking the inner product of both sides of  (\ref{eq:umbcond2}) with $e_3$ yields
\be\label{eq:a0a3} \alpha_0\alpha_3''=y\alpha_3'.\ee
Using (\ref{eq:baralpha}) and (\ref{eq:a0a3}),  we obtain by taking the inner product of both sides of (\ref{eq:umbcond2}) with $\bar{\alpha}$ and $\bar{\alpha}'$, respectively, that
\be\label{eq:ab}z=\alpha_0(\alpha_3')^2\,\,\,\,\mbox{and}\,\,\,\,y=\frac{\alpha_0\varphi'+2\varphi\alpha_0'}{2\varphi}.\ee
Hence (\ref{eq:umbcond}) becomes
\be\label{eq:umbcond3} 2\varphi\alpha_0\alpha''+2\varphi^2e_0-(\alpha_0\varphi'+2\varphi\alpha_0')\alpha'-2\varphi\alpha_0(\alpha_3')^2\bar{\alpha}=0.\ee
  Taking the inner product of both sides of  (\ref{eq:umbcond3}) with $e_3$ yields
  $$\alpha_3''=\left(\frac{\varphi'}{2\varphi}+\frac{\alpha'_0}{\alpha_0}\right)\alpha_3',$$
  which easily implies that 
  \be\label{eq:alpha3}\alpha_3'=c\alpha_0\sqrt{\varphi}\,\,\,\mbox{for some}\,\,\,c\in \R.\ee

We now show that $\alpha_0$ can not be constant on $I$. Assume otherwise, say, that $\alpha_0=a\in \R$.  We may also suppose that $\alpha$ is parametrized by arc-length, i.e., $\varphi=1$. Then $\alpha_3'=ac$ by (\ref{eq:alpha3}), thus $z=a^3c^2$ and $y=0$ by 
(\ref{eq:ab}). Replacing into (\ref{eq:umbcond3}), the $e_0$-component gives $c^2a^4=1$, whereas for $1\leq i\leq 2$ the $e_i$-component then yields $\alpha_i''=(1/a)\alpha_i$. We obtain that $\alpha_i=a_i\exp (s/a)+b_i\exp (-s/a)$ for some $a_i, b_i\in \R$, $1\leq i\leq 2$. Replacing into $1=\alpha_0^2+\alpha_1^2+\alpha_2^2=a^2+\alpha_1^2+\alpha_2^2$ implies that $a_i=0= b_i$ for $1\leq i\leq 2$, i.e., $\alpha_1=0=\alpha_2$, and that $a=\pm 1$. Therefore  $f|_{\tilde U}$ is  totally geodesic, contradicting our assumption. 

 Hence, there must exist an open interval $J\subset I$ such that $\alpha_0'(s)\neq 0$ for all $s\in J$, thus we can reparametrize $\alpha$ on $J$ so that $\alpha_0(s)=s$ for all $s\in J$.  Then the $e_0$-component of (\ref{eq:umbcond3}) gives
$$s\varphi'+2(c^2s^4-1)\varphi^2+2\varphi=0.$$
This is easily seen to be equivalent to   $\varphi^{-1}(s)=c^2s^4+c_2s^2+1$ for some $c_2\in \R$. Hence the right-hand-side of the preceding equation is nowhere vanishing, and we can write $b\varphi^{-1}(s)=s^4+as^2+b$ for $a=c_2/c^2$ and $b=1/c^2$, or equivalently,
 \be\label{eq:varphi}\varphi(s)=\frac{p^2-q}{(s^2-p)^2-q},\,\,\,\,\,p^2>q,\ee for $p=-\frac{a}{2}$ and $q=\frac{a^2}{4}-b$. Equation (\ref{eq:alpha3}) becomes 
\be\label{eq:alpha3b}\alpha_3'(s)=\frac{s}{\sqrt{(s^2-p)^2-q}}.\ee
 Taking the inner product of both sides of  (\ref{eq:umbcond3}) with $e_i$, $1\leq i\leq 2$, and using (\ref{eq:varphi}) and  (\ref{eq:alpha3b}) yields
$$s((s^2-p)^2-q)\alpha_i''+(s^4-p^2+q)\alpha'-s^3\alpha_i=0, \,\,\,1\leq i\leq 2.$$
\begin{lemma}\po\label{le:edo} Let $\alpha_i\colon\,I\to \R$, $1\leq i\leq 2$, be linearly independent solutions of the ODE
\be\label{eq:ODE} s((s^2-p)^2-q)\alpha_i''+(s^4-p^2+q)\alpha'-s^3\alpha_i=0,\,\,\,\,p^2>q,\ee
on an open interval $I\subset \R$ where  $(s^2-p)^2-q>0$. Assume that $s^2+\alpha_1^2+\alpha_2^2=1$ for all $s\in I$. 
Then $(p,q)\in {\cal I}$, $I\subset J_{p,q}$ and there exists $\theta\in \R$ such that
\be\label{eq:sol}
\sqrt{2}(\alpha_1(s),\alpha_2(s))=\left(h(s)+\frac{1-p}{h(s)}\right)(\cos \theta, \sin\theta)\pm \frac{\sqrt{q-(1-p)^2}}{h(s)}(-\sin \theta, \cos\theta),
\ee
where $h(s)=\sqrt{p-s^2+\sqrt{(p-s^2)^2-q}}$.
\end{lemma}
\proof Let $F$ be a primitive of $\beta\colon\, I\to \R$ given by 
$$\beta(s)=\frac{s}{\sqrt{(s^2-p)^2-q}}.$$
Then, it is easily checked that the functions  $$\rho_+:=\exp \circ F\,\,\,\,\mbox{and}\,\,\,\,\rho_-:=\exp\circ ( -F)$$ form a basis of the space of solutions of (\ref{eq:ODE}) on $I$.  Thus, there exist $a_i, b_i\in \R$, $1\leq i\leq 2$, such that
\be\label{eq:alphai}\alpha_i=a_i\rho_++b_i\rho_{-},\,\,\,\,\,\,1\leq i\leq 2.\ee
Replacing into $s^2+\alpha_1^2+\alpha_2^2=1$ gives
\be\label{eq:ABC}s^2+A\exp (2F(s))+ B+ C\exp(-2F(s))=0,\,\,\,\,\mbox{for all}\,\,s\in I,\ee
where  $A=a_1^2+a_2^2$, $B=2(a_1b_1+a_2b_2)-1$ and  $C=b_1^2+b_2^2$. 

Assume that either of the following conditions holds:
$$(i)\,\, q<0; \,\,\,\,\,\,(ii)\,\, q>0\,\,\mbox{and}\,\,p\leq 0;\,\,\,\,\,\,(iii) \,\,\mbox{$q> 0$, $p>0$ and $I$ is not contained in $J_{p,q}$}.$$    
Then, up to a constant, 
$$F(s)=\frac{1}{2}\log\left(s^2-p+\sqrt{(s^2-p)^2-q}\right),$$
hence (\ref{eq:ABC}) gives
$$
A(u+\sqrt{u^2-q})+C(u+\sqrt{u^2-q})^{-1}=-u+E,
$$
where $u=s^2-p$ and $E=-B-p$. This yields
$$2(2A+1)(C-Aq)=-(2A+1)^2q,\,\,\,(C-Aq)E=-(2A+1)Eq\,\,\,\mbox{and}\,\,\,(C-Aq)^2=-E^2q.$$
Since $2A+1>0$, the first and the third of the preceding equations give  $q=0$ if $E=0$, whereas the same conclusion
follows from the first and second equations if $E\neq 0$. This is a contradiction and shows that either $q=0$ or $q> 0$, $p>0$ and $I\subset J_{p,q}$.

    Let us consider first the case $q=0$. Suppose  either that $p<0$ or that $p>0$ and $I$ is not contained in $J_{p,0}=(-\sqrt{p}, \sqrt{p})$.  Then $F(s)=\frac{1}{2}\log (s^2-p)$ and (\ref{eq:ABC}) gives
    $$Au+Cu^{-1}=-u+E,$$
    which implies that $C=0$, $E=0$ and $A=-1$, a contradiction. Thus $p>0$ and $I\subset J_{p,0}$, in which case $F(s)=-\frac{1}{2}\log (p-s^2)$ and (\ref{eq:ABC}) now yields
    $$-Au^{-1}-Cu=-u+E.$$
    This implies that $A=0$, $E=0$ and $C=1$, hence $(a_1, a_2)=(0,0)$, $p=-B=1$ and there exists $\theta\in \R$ such that $(b_1, b_2)=(\cos\theta, \sin\theta)$. Therefore $\alpha_1=\cos\theta\sqrt{1-s^2}$, $\alpha_2=\sin\theta\sqrt{1-s^2}$, and hence the statement is true in this case.
    
    Now suppose that $q> 0$, $p>0$ and $I\subset J_{p,q}$. Then 
    \be\label{eq:F}F(s)=-\frac{1}{2}\log\left(\sqrt{(p-s^2)-q}-s^2+p\right)=-\log h(s) ,\ee
  where $h(s)$ is as in the statement. We obtain from (\ref{eq:ABC}) that
 $$
A(\sqrt{u^2-q}-u)^{-1}+C(\sqrt{u^2-q}-u)=-u+E,
$$
with $u=s^2-p$ and $A,C,E$ as before. This is equivalent to
$$2(2C-1)(A-qC)=-q(2C-1)^2,\,\,\,\,(A-qC)E=-q(2C-1)E\,\,\,\mbox{and}\,\,\,(A-qC)^2=-qE^2,$$
and hence to
$$E=0, \,\,\,A=\frac{q}{2}\,\,\,\mbox{and}\,\,\,C=\frac{1}{2}.$$
This gives
\be\label{eq:aibi}
a_1b_1+a_2b_2=\frac{1}{2}(1-p),\,\,\,\,\,a_1^2+a_2^2=\frac{q}{2}\,\,\,\,\mbox{and}\,\,\,\,b_1^2+b_2^2=\frac{1}{2}.
\ee  
By the last  equation in (\ref{eq:aibi}),  there exists $\theta\in \R$ such that $(b_1, b_2)=\frac{\sqrt{2}}{2}(\cos\theta, \sin\theta)$.
Set $u:=(\cos\theta, \sin\theta)$ and $v:=(-\sin\theta, \cos\theta)$. Then the first  equation can be written as 
\be\label{eq:u}\<(a_1,a_2), u\>=\frac{\sqrt{2}}{2}(1-p).\ee Using this and the second equation we obtain
 $$\frac{q}{2}=a_1^2+a_2^2=\frac{1}{2}(p-1)^2+\<(a_1,a_2), v\>^2,$$
 hence
 $$\<(a_1,a_2), v\>^2=\frac{1}{2}(q-(p-1)^2).$$
 In particular, this shows that  $q\geq (p-1)^2$, thus  $(p,q)\in {\cal I}$. Moreover, together with (\ref{eq:u}) it implies  that 
 $$a_1=\frac{\sqrt{2}}{2}((1-p)\cos \theta\mp \sqrt{q-(1-p)^2}\sin \theta),\,\,\,\,a_2=\frac{\sqrt{2}}{2}((1-p)\cos \theta\pm \sqrt{q-(1-p)^2}\sin \theta),$$
and (\ref{eq:sol}) follows.\vspace{1ex}\qed
 
 We obtain from Lemma \ref{le:edo} and (\ref{eq:alpha3b}) that $\alpha(s)$ is given by
 $$\left(s,\frac{\sqrt{2}}{2}\left(h(s)+\frac{1-p}{h(s)}\right)(\cos \theta, \sin\theta)\pm \frac{\sqrt{2}}{2}\frac{\sqrt{q-(1-p)^2}}{h(s)}(-\sin \theta, \cos\theta),h\pm \log h(s)\right),$$
 and $f|_{\tilde U}$ can be parametrized by $Y\colon\, \Sf^{m-1}\times I\to \Sf^{m+1}\times \R$ given by $$Y(X,s)=(sX,\alpha_1(s), \alpha_2(s),\alpha_3(s)).$$ Let $A$ be the linear isometry of $\R^{m+2}\supset \Sf^{m+1}\times \R$ defined by
 $$Ae_m=\cos\theta e_m +\sin \theta e_{m+1}, \,\,\,\,Ae_{m+1}=\mp \sin\theta e_m +\cos \theta e_{m+1},$$
 $Ae_i=e_i$ for $i\in \{1,\ldots, m\}$ and $Ae_{m+2}=\pm e_{m+2}$. Then $A^{-1}Y(X,s)-he_{m+2}=Y_{p,q}(X,s)$.\vspace{1ex}
 
 It remains to prove assertion $(v)$ in the direct statement. This is equivalent to  showing that $Y_{p,q}$ and $Y_{p',q'}$ do  not parametrize congruent submanifolds for distinct pairs $(p,q)$ and $(p', q')$ in ${\cal I}$. After reparametrizing the curve $\alpha=\alpha_{p,q}$ by arc-length, 
 the metric induced by $Y_{p,q}$ is a warped product metric $ds^2+\rho^2(s)d\sigma$ on $I\times \Sf^{m-1}$, where $d\sigma$ is the
 standard metric on $\Sf^{m-1}$ and the warping function $\rho=\rho_{p,q}$ is the inverse of the arc-length function $$S_{p,q}(t)=\int_0^t\|\alpha_{p,q}'(\tau)\|d\tau= \int_0^t\sqrt{\varphi_{p,q}(\tau)}d\tau,$$ with $\varphi_{p,q}$ given by (\ref{eq:varphi}).  If $Y_{p,q}$ and $Y_{p',q'}$ parametrize congruent submanifolds, then the induced metrics, and hence the corresponding warping functions, must  coincide. It follows that $\varphi_{p,q}=\varphi_{p',q'}$, which easily  implies that $(p,q)=(p', q')$.\vspace{2ex}\qed
 
\noindent {\em Acknowledgement:} A first draft of Lemma \ref{le:redcod} was derived in a conversation of the second author with M. Dajczer. We thank him for allowing us to include it here.

\vspace*{-1ex} {\renewcommand{\baselinestretch}{1}
\hspace*{-25ex}\begin{tabbing}
\indent  \=  Universidade Federal de S\~ao Carlos \\
\indent  \= Via Washington Luiz km 235 \\
\> 13565-905 -- S\~ao Carlos -- Brazil \\
\> e-mail: bruno@dm.ufscar.br \\
\hspace*{10ex} tojeiro@dm.ufscar.br
\end{tabbing}}

\end{document}